\documentclass[aip,reprint]{revtex4-1}

\usepackage{amsfonts,amsmath,amssymb}
\usepackage{mathtools}
\usepackage{commath}
\usepackage[none]{hyphenat}
\usepackage{fancyhdr}
\usepackage{graphicx}
\usepackage{subfigure}
\usepackage{float}
\usepackage{comment}
\usepackage{bm}
\usepackage{xcolor}
\usepackage{enumerate}
\usepackage{enumitem}
\usepackage{xspace}
\usepackage{algorithm}
\usepackage{algpseudocode}
\usepackage{currvita}
\usepackage{ragged2e}
\usepackage{amsthm}
\usepackage{amssymb}
\usepackage{hyperref}
\usepackage{tikz}
\usetikzlibrary{decorations.pathmorphing,patterns}
\usetikzlibrary{calc,patterns,decorations.markings}
\usetikzlibrary{positioning}
\usetikzlibrary{matrix,positioning,decorations.pathreplacing}
\usetikzlibrary{matrix, fit}
\usepackage[english]{babel}
\newtheorem{theorem}{Theorem}
\newtheorem{proposition}[theorem]{Proposition}
\newtheorem{assumption}[theorem]{Assumption}

\draft 

\begin{document}


\title{The spatiotemporal coupling in delay-coordinates dynamic mode decomposition} 



\author{Emil Bronstein}
\email[]{emilbr@campus.technion.ac.il}
\affiliation{Faculty of Mechanical Engineering, Technion -- Israel Institute of Technology, Haifa, Israel}

\author{Aviad Wiegner}
\affiliation{Viterbi Faculty of Electrical and Computer Engineering, Technion -- Israel Institute of Technology, Haifa, Israel}

\author{Doron Shilo}
\affiliation{Faculty of Mechanical Engineering, Technion -- Israel Institute of Technology, Haifa, Israel}

\author{Ronen Talmon}
\affiliation{Viterbi Faculty of Electrical and Computer Engineering, Technion -- Israel Institute of Technology, Haifa, Israel}


\date{\today}


\begin{abstract}
Dynamic mode decomposition (DMD) is a leading tool for equation-free analysis of high-dimensional dynamical systems from observations.
In this work, we focus on a combination of delay-coordinates embedding and DMD, i.e., \textit{delay-coordinates DMD}, which accommodates the analysis of a broad family of observations. An important utility of DMD is the compact and reduced-order spectral representation of observations in terms of the DMD eigenvalues and modes, where the temporal information is separated from the spatial information.
From a spatiotemporal viewpoint, we show that when DMD is applied to delay-coordinates embedding, temporal information is intertwined with spatial information, inducing a particular spectral structure on the DMD components. We formulate and analyze this structure, which we term the \textit{spatiotemporal coupling in delay-coordinates DMD}.
Based on this spatiotemporal coupling, we propose a new method for DMD components selection. When using delay-coordinates DMD that comprises redundant modes, this selection is an essential step for obtaining a compact and reduced-order representation of the observations.
We demonstrate our method on noisy simulated signals and various dynamical systems and show superior component selection compared to a commonly-used method that relies on the amplitudes of the modes.
\end{abstract}


\maketitle 

\begin{quotation}
Dynamical systems are abundant in many fields of science and engineering. As dynamical systems are often high-dimensional and complex, their characterization from observations is a coveted goal. A key task in accomplishing this goal is finding a reduced-order representation of a system, where dynamic mode decomposition (DMD) is a prominent tool for this purpose. In recent years, the combination of delay-coordinates embedding and DMD (i.e., delay-coordinates DMD) has been shown to be highly useful in the characterization of dynamical systems, even when these systems are highly nonlinear or chaotic. This combination gives rise to a specific spectral structure that does not exist in ordinary DMD and has not been studied so far. In this work, we formulate and analyze this structure, where we show that the spectral components of delay-coordinates DMD exhibit a spatiotemporal coupling. This coupling suggests that the representations obtained by delay-coordinates DMD can be further reduced. By considering this coupling, we propose not only to decouple the spatial information from the temporal information but also to exploit it to construct an improved reduced-order representation in an unsupervised fashion. We demonstrate our approach on several dynamical systems that include noisy observations of undamped and damped mechanical oscillators.
\end{quotation}



\section{Introduction} \label{sec:introduction}

Time-series analysis, modeling, and prediction are ubiquitous tasks in applied sciences. When the time-series stem from ergodic dynamical systems, learning their phase space in a nonparametric fashion from sufficiently long intervals of observations is possible and has become an active field of research in recent years. 
This so-called data-driven approach seeks to obtain meaningful, physics-related knowledge, in an \textit{equation-free} manner, shifting the focus from equations-based descriptions to observations-based analysis \cite{kevrekidis2004equation,kevrekidis2009equation,kevrekidis2003equation,sirisup2005equation,theodoropoulos2000coarse,gear2002coarse}.
Existing data-driven methods can be largely divided into two categories \cite{budivsic2012applied}.
The first is based on a state-space representation, and the focus is typically on finding a map that propagates a state at a given time to a state at a future time. Arguably, the classical approach is to approximate the nonlinear dynamics as a collection of locally linear systems on tangent spaces near attractors \cite{farmer1987predicting,sauer1994time,kugiumtzis1998regularized}.
Another, more recent class of methods in this category, constructs reduced models of the state space in Euclidean spaces, and then, finds the propagation rules in the reduced spaces \cite{kevrekidis2004equation,talmon2013empirical,dsilva2016data,crosskey2017atlas}.

The second category is an operator-theoretic approach that considers observables of the system.
In this approach, many recent methods are based on the Koopman operator\cite{koopman1931hamiltonian}, which is a linear, infinite-dimensional operator that operates on observables in a Hilbert space, and propagates them \emph{linearly} in time \cite{mezic2013analysis,kutz2016dynamic}.
While the ability to represent nonlinear dynamics in a linear fashion is perhaps the most notable property of using the Koopman operator, it has several other remarkable attributes. 
In this line of work\cite{mezic2004comparison,mezic2005spectral,mezic2013analysis,arbabi2017ergodic,budivsic2012applied}, it was shown that the Koopman operator has the ability of capturing the  dynamics of linear or nonlinear dynamical systems through its spectral components.
More concretely, the time evolution of a dynamical system can be decomposed into spatial patterns, which are often referred to as the Koopman modes and derived from the eigenfunctions of the Koopman operator, and temporal patterns, which are derived from the eigenvalues of the Koopman operator. Therefore, the spectral analysis of the Koopman operator is of high importance for obtaining an informative description of dynamical systems \cite{giannakis2019data,budivsic2012applied,mezic2004comparison,mezic2005spectral,pan2021sparsity,das2020koopman}.

Using the Koopman operator for analysis poses an important tradeoff. On the one hand, it facilitates a linear description of the dynamics and a useful decoupling of the spatial and temporal patterns via its spectral representation. On the other hand, the Koopman operator is infinite dimensional. Therefore, finite-dimensional approximations are necessary for practical purposes. Perhaps the most common technique for such a finite approximation is the dynamic mode decomposition (DMD), introduced by Schmid and Sesterhenn\cite{schmid2009dynamic,schmid2010dynamic}, where the primary goal is to approximate the Koopman eigenvalues and modes in finite spaces based on a finite set of observations \cite{rowley2009spectral}.
The DMD eigenvalues and modes give rise to a \emph{discrete} spatiotemporal representation of the dynamical system, facilitating a \emph{data-driven} and \emph{equation-free} analysis. 
Following its appearance, several variants for enhancing the capabilities of DMD have been introduced, e.g., the extended DMD (EDMD)\cite{williams2015data}, where projections on finite-space dictionaries were used to improve the spatiotemporal representation, least-squares \cite{chen2012variants,dawson2016characterizing} and sparsity promoting \cite{jovanovic2014sparsity,sayadi2015parametrized} techniques, as well as approaches for the analysis of compressed data \cite{brunton2015compressed,brunton2014compressive}.
Additional extensions of DMD have been introduced, such as multiresolution DMD (mrDMD)\cite{kutz2016multiresolution}, which effectively uncovers multiscale structures in the data, and DMD with control (DMDc) \cite{proctor2016dynamic}, which extracts low-order models of high-dimensional systems that require control.

Indeed, in recent years, DMD has been shown to be a powerful tool, demonstrating remarkable capabilities in a broad range of fields \cite{brunton2017chaos,brunton2019methods,mann2016dynamic,elmore2021learning,brunton2016extracting,mohan2018data}, and particularly, in fluid dynamics \cite{rowley2009spectral,duke2012error,muld2012flow,seena2011dynamic,tu2014spectral,bagheri2013koopman,rowley2017model}.
Nevertheless, despite its popularity and success, DMD has some notable limitations. For example, a straightforward application of DMD does not allow for the signal reconstruction of a standing wave  \cite{kutz2016dynamic,tu2014on}.
In addition, ordinary DMD cannot handle cases when the number of linearly-independent DMD modes is smaller than the number of the system's oscillation frequencies \cite{le2017higher}. 
Another notable limitation of DMD concerns one-dimensional signals, as its application to such signals results in
a degenerate (scalar) representation.

Combining delay-coordinates embedding and DMD mitigates these limitations\cite{kutz2016dynamic,tu2014on,le2017higher}.
The delay-coordinates embedding (also referred to as time-delay embedding, or, simply, delay-coordinates) is a method for augmenting past observations to the present observation. This approach dates back to 1981 with the formulation of Takens' embedding theorem\cite{takens1981detecting}, according to which an attractor of a system can be reconstructed up to a diffeomorphism using delay embedding.

Various studies that associate delay-coordinates and DMD-based methods have been presented in recent years. For example, Le Clainche and Vega\cite{le2017higher} presented the higher order DMD (HODMD), a global linear method that is capable of uncovering a large number of frequencies of periodic and quasiperodic dynamical systems based on limited and noisy input data. Brunton \textit{et al}. introduced the HAVOK analysis in terms of a Hankel matrix, which successfully represents highly nonlinear and chaotic systems using a linear model and intermittent forcing \cite{brunton2017chaos}. Pan and Duraisamy provided the minimal required augmentation number for a perfect recovery of dynamical systems based on their Fourier spectrum \cite{pan2020structure}. 
We note that combinations of delay-coordinates and Koopman operator-based methods have also been investigated\cite{arbabi2017ergodic,das2019delay,kamb2020time}.

In this paper, we show that despite its prevalence, application of DMD to \textit{augmented} data (i.e., delay-coordinates DMD) gives rise to a particular spectral structure that does not exist in ordinary DMD applications and, to the best of our knowledge, has not been studied in the existing literature. This structure comprises augmented DMD modes that embody \emph{temporal} information, which is entangled with \emph{spatial} information. We formulate and analyze this entanglement, and term it as the \emph{spatiotemporal coupling in delay-coordinates DMD}.

Based on this spatiotemporal coupling, we propose a new approach for obtaining compact and reduced-order representations of dynamical systems from observations, where, similarly to ordinary DMD, the spatial and temporal patterns are \emph{decoupled}. 
Our approach includes solving an inherent challenge of delay-coordinates DMD, where the number of \emph{augmented} DMD components is often larger than the number of intrinsic modes of the dynamical system. In such cases, the augmented DMD components can be divided into two subsets: those that describe the dynamical system, which we term \emph{true}, and those that are a mere artifact of the augmentation, unrelated to the system, which we term \emph{spurious}.

By relying on the difference in their spatiotemporal coupling, we distinguish between the true and spurious DMD components, and present a method for selecting the true components. In contrast to ordinary DMD, our method is based on delay-coordinates DMD, and is therefore capable of analyzing and representing a broader family of observations and signals. We demonstrate the spatiotemporal coupling and the effectiveness of our approach on various simulated dynamical systems.

This paper is organized as follows. Section \ref{sec:related_work} briefly describes existing approaches for mode selection. Section \ref{sec:prob_formulation} presents the problem formulation. In Section \ref{sec:sdcdmd}, we introduce the spatiotemporal coupling in delay-coordinates DMD in detail, and reveal the specific relations within augmented DMD modes.
Then, in Section \ref{sec:proposed_method}, we propose a new method for decoupling the spatial and temporal information, leading to the identification of the true augmented DMD components and to a compact,  reduced-order and informative representation of the observations.
For the purpose of illustration, in Section \ref{sec:example}, we present the application of our method to a two-mode sine signal.
Finally, in Section \ref{sec:results}, our method is demonstrated on various dynamical systems, outperforming the common method that relies on the amplitudes of the modes at low SNR values.


\section{Related work on mode selection}\label{sec:related_work}

The identification of the dominant DMD components required for optimal, reduced-order representations of dynamical systems (which we term \emph{true}) has been studied in wider contexts, beyond delay-coordinates DMD.
Criteria and methods for such a purpose are usually referred to as \emph{mode selection}, where we use the broader term \emph{DMD components selection}.
For instance, Rowley \textit{et al.} ordered the DMD modes by their norms (amplitudes)\cite{rowley2009spectral} -- a method we term the \emph{maximal amplitudes method}.
The norms of the modes can be further weighted by the magnitudes of the corresponding eigenvalues to account for modes that have large norms yet decay rapidly, as suggested by Tu \textit{et al.} \cite{tu2014on}.
Schmid \textit{et al.} used a projection of the data sequence on the identified modes, whose coefficients indicate on the significance of the modes \cite{schmid2012decomposition}.
Jovanovi\'c \textit{et al.} introduced a sparsity-promoting approach using an addition of a penalty term of the DMD amplitudes, which regularizes the least-square deviation between the linear combination of DMD modes and the snapshots matrix\cite{jovanovic2014sparsity}.
Tissot \textit{et al.} proposed an energetic criterion, where the amplitude of a DMD mode is weighted by its corresponding temporal coefficient \cite{tissot2014model}.
Sayadi \textit{et al.} proposed a parameterized approach, which utilizes the sparsity promoting DMD, and then reconstructs the modes' amplitudes using a time-dependent coefficient, giving a notion of their significance \cite{sayadi2015parametrized}. 
Another approach was proposed by Kou and Zhang, which considered the initial conditions and temporal evolution of DMD modes, ordering them according to the integrals of their corresponding time coefficients \cite{kou2017improved}.

Our work differs from the above studies on mode selection,
as it provides, for the first time to the best of our knowledge, a mode selection framework in the specific context of delay-coordinates DMD.
Seemingly, the problem of mode selection in delay-coordinates DMD is more challenging due to the existence of spurious DMD modes, as well as the higher dimensions of augmented DMD modes compared to the observations. Nevertheless, we show that despite the additional challenges, delay-coordinates DMD also constitutes a remedy. Specifically, we show that the spatiotemporal coupling in delay-coordinates DMD bears information that facilitates a new method for mode selection.


\section{Problem Formulation} \label{sec:prob_formulation}
Consider a dynamical system
\begin{equation}
    \bm{x}_{k+1}
    =
    \bm{f}(\bm{x}_k),  \quad k\in \mathbb{Z},
    \label{eq:dynamical_system}
\end{equation}
that evolves in discrete time $k$ on a manifold $\mathcal{M}  \subset \mathbb{R}^n$, where $\bm{x}_k \in \mathcal{M}$ is a state vector, and $\bm{f}: \mathcal{M} \mapsto \mathcal{M}$. The function $\bm{f}$ is unknown and could be deterministic or stochastic, e.g., due to the presence of noise.
Assume that this discrete time formulation arises from a continuous-time dynamical system $\dot{\bm{x}} = \mathcal{\bm{A}} \bm{x}$. Specifically, suppose that the state is sampled at a fixed sampling rate $\omega_s=2\pi / \Delta t$, and $m+1$ discrete samples $\bm{x}_k=\bm{x}(k \Delta t)$ of the state are collected, where $ k = 0,1,\dots,m$.
In this work, following common practice \cite{le2017higher,pan2020structure,kamb2020time}, we assume that $m>n$.

Our goal is to analyze the dynamics of system \eqref{eq:dynamical_system} based on the finite set of observations $\{\bm{x}_k\}_{k=0}^{m}$.
For this purpose, we use the DMD approach \cite{schmid2009dynamic,budivsic2012applied,tu2014on}, which provides a spatiotemporal representation of the observations, and is described next.
First, a linear approximation of the discrete-time system in \eqref{eq:dynamical_system} is employed:
\begin{equation}
    \bm{x}_{k+1} = \bm{A} \bm{x}_k, \quad \bm{A} \in \mathbb{R}^{n \times n}, \quad \bm{x}_k \in \mathbb{R}^n.
    \label{eq:A}
\end{equation}
To find $\bm{A}$, the observations $\{\bm{x}_k\}_{k=0}^m$ are arranged into two observation matrices
\begin{equation}
\begin{split}
    \bm{X}
    & =
    \begin{bmatrix}
    | & | &  & | \\
    \bm{x}_0 & \bm{x}_1 & \cdots & \bm{x}_{m-1} \\
    | & | &  & |
    \end{bmatrix}  \in \mathbb{R}^{n \times m}, \\
    \bm{X}'
    & =
    \begin{bmatrix}
    | & | &  & | \\
    \bm{x}_1 & \bm{x}_2 & \cdots & \bm{x}_m \\
    | & | &  & |
    \end{bmatrix}  \in \mathbb{R}^{n \times m}.
    \end{split}
    \label{eq:obs_mat}
\end{equation}

Then, following the exact DMD \cite{tu2014on,kutz2016dynamic}, the singular value decomposition (SVD) of $\bm{X}$ is computed as
\begin{equation}
    \bm{X} = \bm{U} \bm{\Sigma} \bm{V}^*,
\end{equation}
where $\bm{U} \in \mathbb{C}^{n \times r}$, $\bm{\Sigma} \in \mathbb{C}^{r \times r}$, $\bm{V} \in \mathbb{C}^{m \times r}$, $( \cdot )^*$ is the complex conjugate transpose, and $r = \text{rank}(\bm{X})$. 
Define $\Tilde{\bm{A}} \in \mathbb{C}^{r \times r}$  by
\begin{equation}
    \Tilde{\bm{A}} = \bm{U}^* \bm{X}' \bm{V} \bm{\Sigma}^{-1}.
    \label{eq:DMD_eigenvalues}
\end{equation}
The so-called DMD eigenvalues, $\lambda_i$, are the eigenvalues of $\Tilde{\bm{A}}$  that satisfy
\begin{equation}
    \Tilde{\bm{A}} \bm{v}_i = \lambda_i \bm{v}_i, \quad i=1, \dots, r,
\end{equation}
where $\bm{v}_i$ are the corresponding eigenvectors.
The so-called exact DMD modes, $\bm{\phi}_i$, are given by\cite{tu2014on,kutz2016dynamic}
\begin{equation}
    \bm{\phi}_i = \frac{1}{\lambda_i} \bm{X}' \bm{V} \bm{\Sigma}^{-1} \bm{v}_i, \quad i=1, \dots, r.
    \label{eq:exact_modes}
\end{equation}
Tu \textit{et al}. defined the exact DMD modes, $\bm{\phi}_i$, and showed that $\left\{ \lambda_i, \bm{\phi}_i \right\} $ are the eigenvalue-eigenvector pairs of $\bm{A}$ \cite{tu2014on}. Moreover, the authors distinguished between the exact modes $\bm{\phi}_i$ and the DMD modes obtained by the standard DMD \cite{schmid2010dynamic} (termed projected DMD modes). It was noted that the exact and projected modes have the tendency to converge when $\bm{X}$ and $\bm{X}'$ have the same column spaces \cite{kutz2016dynamic}. 

The observations $\bm{x}_k$ can be expressed as
\begin{equation}
    \bm{x}_k = \sum_{i=1}^r \lambda_i^k \bm{\phi}_i \eta_{0,i},
    \label{eq:xk}
\end{equation}
where $\eta_{0,i}$ is the $i$th entry of $\bm{\eta}_0 = \bm{\Phi}^* \bm{x}_0$, and $\bm{\Phi}$ is the matrix whose columns are the $r$ leading DMD modes $\bm{\phi}_i$, satisfying $\bm{\Phi}^*\bm{\Phi}=\bm{I}$.
Alternatively, \eqref{eq:xk} can be written in a matrix form as
\begin{equation}
\begin{split}
    \bm{x}_k & = \bm{\Phi} \bm{c}_k,
    \\
    \bm{\Phi} & = 
    \begin{bmatrix}
    | & | &  & | \\
    \bm{\phi}_1 & \bm{\phi}_2 & \cdots & \bm{\phi}_{r} \\
    | & | &  & |
    \end{bmatrix}
    \in \mathbb{C}^{n \times r},
    \\
    \bm{c}_k & = 
    \begin{bmatrix}
    \lambda_1^k \eta_{0,1} \\
    \lambda_2^k \eta_{0,2} \\
    \vdots \\
    \lambda_{r}^k \eta_{0,r} \\
    \end{bmatrix}
    \in \mathbb{C}^{r}.
    \end{split}
    \label{eq:xkmatrix_unaugmented}
\end{equation}

The continuous-time counterpart of the decomposition in \eqref{eq:xk}, which corresponds to the dynamical system $\dot{\bm{x}} = \mathcal{\bm{A}} \bm{x}$, is given by \cite{kutz2016dynamic}
\begin{equation}
    \bm{x}(t) = \sum_{i=1}^r \exp(\mu_i t) \bm{\psi}_i \nu_{0,i},
    \label{eq:xt}
\end{equation}
where $\mu_i$ and $\bm{\psi}_i$ are the eigenvalues and eigenvectors of the propagating operator $\mathcal{A}$, $\nu_{0,i} = \langle \bm{\psi}_i, \bm{x}(0) \rangle =  \bm{\psi}_i^* \bm{x}(0)$, and 
\begin{equation}
    \lambda_i = \exp(\mu_i \Delta t).
    \label{eq:lambdamu}
\end{equation}

In case of a damped oscillating system, when the system comprises natural frequencies $\omega_i$ and damping ratios $\zeta_i$, according to \eqref{eq:lambdamu}, we have
\begin{equation}
    \mu_i = -\zeta_i \omega_i \pm j\omega_i \sqrt{1-\zeta^2},
    \quad
    j = \sqrt{-1},
    \label{eq:mu_zeta_omega}
\end{equation}
from which we obtain
\begin{equation}
    \omega_i =  \frac{ \left| \log(\lambda_i) \right| }{\Delta t},
    \quad
    \zeta_i = -\frac{\Re( \log(\lambda_i) )}{\omega_i\Delta t }.
    \label{eq:omegafromlambda}
\end{equation}
In case the system is undamped, $\zeta=0$ can be substituted into \eqref{eq:mu_zeta_omega} and \eqref{eq:omegafromlambda}.

Eq. \eqref{eq:xk} constitutes a decomposition of the system in \eqref{eq:dynamical_system} into its dynamic modes, such that $\lambda_i$ and $\bm{\phi}_i$ hold the temporal and spatial information, respectively.
Yet, as described in Section \ref{sec:introduction}, obtaining this decomposition might be hindered for a variety of reasons, such as standing waves, one- or low-dimensional signals\cite{kutz2016dynamic,tu2014on},
and observation noise (e.g., noise in the measurement equipment).
Our goal is to obtain a compact and reduced-order representation of the system similar to Eq. \eqref{eq:xk}, where the challenges mentioned above are present. To accomplish this goal, we use \textit{delay-coordinates DMD} as detailed below.


\section{Spatiotemporal coupling in delay-coordinates DMD} \label{sec:sdcdmd}

In this section, we show that when DMD is applied to delay-coordinates embedding (constituting the delay-coordinates DMD), the strict separation of temporal and spatial information in representations obtained by DMD as in \eqref{eq:xk} is violated.
In other words, the separation where the eigenvalues $\lambda_i$ bear the temporal information about the dynamical system and the modes $\bm{\phi}_i$ represent the spatial information no longer exists in delay-coordinates DMD.
Then, we show that this seemingly limiting spatial and temporal information entanglement could be harnessed toward the selection of the DMD components required for accurate and reduced-order characterization of the dynamical system.
Specifically, we formulate and analyze the induced spatiotemporal structure of the eigenvalues and modes that arise from delay-coordinates DMD, which we term \textit{augmented} DMD components. For simplicity, we divide the exposition into two stages.
First, an augmentation of one sample is considered in Subsection \ref{subsec:s=1}, followed by a generalization to augmentations of several samples, which is presented in Subsection \ref{subsec:s>1}.

We begin by formulating the augmentation, i.e., applying delay-coordinates embedding to the observation matrices $\bm{X}$ and $\bm{X}'$ in Eq. \eqref{eq:obs_mat}.
Let $\hat{\bm{X}}$ and $\hat{\bm{X}'}$ denote the augmented observation matrices, defined by 
\begin{equation}
\begin{split}
    \hat{\bm{X}}
    & =
    \begin{bmatrix}
    | & | &  & | \\
    \hat{\bm{x}}_0 & \hat{\bm{x}}_1 & \cdots & \hat{\bm{x}}_{m-s-1} \\
    | & | &  & |
    \end{bmatrix}
    \\
    & =
    \begin{bmatrix}
    \bm{x}_0 & \bm{x}_1 & \cdots & \bm{x}_{m-s-1} \\
    \bm{x}_1 & \bm{x}_2 & \cdots & \bm{x}_{m-s} \\
    \vdots & \vdots & \vdots & \vdots \\
    \bm{x}_s & \bm{x}_{s+1} & \cdots & \bm{x}_{m-1}
    \end{bmatrix}
    \in \mathbb{R}^{[n(s+1)] \times [m-s]},
    \\
    \hat{\bm{X}'}
    & =
    \begin{bmatrix}
    | & | &  & | \\
    \hat{\bm{x}}_1 & \hat{\bm{x}}_2 & \cdots & \hat{\bm{x}}_{m-s} \\
    | & | &  & |
    \end{bmatrix}
    \\
    & =
    \begin{bmatrix}
    \bm{x}_1 & \bm{x}_2 & \cdots & \bm{x}_{m-s} \\
    \bm{x}_2 & \bm{x}_3 & \cdots & \bm{x}_{m-s+1} \\
    \vdots & \vdots & \vdots & \vdots \\
    \bm{x}_{s+1} & \bm{x}_{s+2} & \cdots & \bm{x}_{m}
    \end{bmatrix}
    \in \mathbb{R}^{[n(s+1)] \times [m-s]},
    \label{eq:aug_obs_mat}
\end{split}
\end{equation}
where $s < m$, $s \in \mathbb{N}$ is termed the \textit{augmentation number}. 
In the remainder of the paper, we use hats to denote either augmented or augmented-related terms. Note that $\hat{\bm{X}}$ and $\hat{\bm{X}'}$ are Hankel matrices, which are typical in delay-coordinates DMD \cite{brunton2017chaos,kutz2016dynamic,kamb2020time}.

Same as in the exact DMD \cite{tu2014on,kutz2016dynamic}, $\hat{\bm{X}}$ and $\hat{\bm{X}'}$ are related through
\begin{equation}
\hat{\bm{X}'} = \hat{\bm{A}} \hat{\bm{X}},
    \quad
    \hat{\bm{A}} \in \mathbb{R}^{[n(s+1)] \times [n(s+1)]},
\end{equation}
or in vector form,
\begin{equation}
\begin{split}
    \hat{\bm{x}}_{k+1} & = \hat{\bm{A}} \hat{\bm{x}}_k,
    \quad
    \hat{\bm{x}}_k, \hat{\bm{x}}_{k+1} \in \mathbb{R}^{n(s+1)},
    \\
    k & = 0, \dots, m-s-1.
\end{split}
    \label{eq:Aaug}
\end{equation}

Application of the exact DMD to $\hat{\bm{X}}$ and $\hat{\bm{X}}'$ results in $\hat{r}$ DMD eigenvalue-mode pairs denoted by $\{ \hat{\lambda}_l, \hat{\bm{\phi}}_l \}_{l=1}^{\hat{r}}$, where $\hat{r}=\text{rank}(\hat{\bm{X}})$.
Then, similarly to Eq. \eqref{eq:xk}, the augmented observations $\hat{\bm{x}}_k$ can be represented as
\begin{equation}
    \hat{\bm{x}}_k = \sum_{l=1}^{\hat{r}} \hat{\lambda}_l^k \hat{\bm{\phi}}_l \hat{\eta}_{0,l},
    \quad
    k=0, \dots, m-s-1,
    \label{eq:xkbasis}
\end{equation}
where $\hat{\eta}_{0,l}$ is the $l$th entry of $\hat{\bm{\eta}}_0 = \hat{\bm{\Phi}}^* \hat{\bm{x}}_0$, and $\hat{\bm{\Phi}}$ is the matrix whose columns are the augmented DMD modes $\hat{\bm{\phi}}_l$. Eq. \eqref{eq:xkbasis} can be written in a matrix form as
\begin{equation}
\begin{split}
    \hat{\bm{x}}_k & = \hat{\bm{\Phi}} \hat{\bm{c}}_k,
    \\
    \hat{\bm{\Phi}} & = 
    \begin{bmatrix}
    | & | &  & | \\
    \hat{\bm{\phi}}_1 & \hat{\bm{\phi}}_2 & \cdots & \hat{\bm{\phi}}_{\hat{r}} \\
    | & | &  & |
    \end{bmatrix}
    \in \mathbb{C}^{[n(s+1)] \times \hat{r}},
    \\
    \hat{\bm{c}}_k & = 
    \begin{bmatrix}
    \hat{\lambda}_1^k \eta_{0,1} \\
    \hat{\lambda}_2^k \eta_{0,2} \\
    \vdots \\
    \hat{\lambda}_{\hat{r}}^k \eta_{0,\hat{r}} \\
    \end{bmatrix}
    \in \mathbb{C}^{\hat{r}}.
\end{split}
    \label{eq:xkmatrix}
\end{equation}
In order to obtain a sufficient amount of augmented DMD components for a full description of the dynamical system, $\hat{r}$ must be at least as large as the number of DMD components required for the description of the unaugmented (original) system in \eqref{eq:xk}; i.e., $\hat{r} \geq r$.
Recalling that $\hat{r}  = \text{rank}(\hat{\bm{X}}) \leq \min \left(n(s+1),m-s \right)$, fulfillment of  $\hat{r} \geq r$ depends on the choice of $s$.
Concretely, choosing $s$ such that $n(s+1) \leq m-s$ leads to the range
\begin{equation}
    0 \leq s \leq \frac{m-n}{n+1}.
    \label{eq:s_range_small}
\end{equation}
Conversely, the choice $n(s+1) \geq m-s$ imposes $r \leq \hat{r} \leq m-s$, which leads to
\begin{equation}
    \frac{m-n}{n+1} \leq s \leq m-r.
    \label{eq:s_range}
\end{equation}
In practice, as $r= \text{rank}(\bm{X})$ is unknown prior to the application of DMD, one can choose $s \leq m-n$ as the upper bound in \eqref{eq:s_range}.

By \eqref{eq:Aaug} and by recalling that $\hat{\bm{A}} \hat{\bm{\phi}}_l = \hat{\lambda}_l \hat{\bm{\phi}}_l$, the expansion of the augmented observations vector $\hat{\bm{x}}_{k+1}$ is given by
\begin{equation}
    \hat{\bm{x}}_{k+1} = \hat{\bm{A}} \hat{\bm{x}}_k =
    \sum_{l=1}^{\hat{r}} \hat{\lambda}_l^k \hat{\lambda}_l \hat{\bm{\phi}}_l \hat{\eta}_{0,l},
    \label{eq:xk1basis}
\end{equation}
which, using \eqref{eq:xkmatrix}, can be recast in matrix form as
\begin{equation}
\begin{split}
    \hat{\bm{x}}_{k+1} & = \hat{\bm{\Phi}} \hat{\bm{\Lambda}} \hat{\bm{c}}_k,
    \\
    \hat{\bm{\Lambda}} & = \text{diag}[\hat{\lambda}_1, \hat{\lambda}_2, \dots, \hat{\lambda}_{\hat{r}}] \in \mathbb{C}^{\hat{r} \times \hat{r}}.
\end{split}
    \label{eq:xk1matrix}
\end{equation}

Eqs. \eqref{eq:xkbasis} and \eqref{eq:xk1basis} show that the temporal propagation of $\hat{\bm{x}}_k$ can be expressed using the augmented DMD components. Specifically, multiplication of the augmented modes $\hat{\bm{\phi}}_l$ by their corresponding augmented eigenvalues $\hat{\lambda}_l$ propagates $\hat{\bm{x}}_k$ one time step to $\hat{\bm{x}}_{k+1}$.
Similarly, $\hat{\bm{x}}_k$ can be propagated $q$ steps to $\hat{\bm{x}}_{k+q}$ by
\begin{equation}
    \hat{\bm{x}}_{k+q} =
    \hat{\bm{A}}^q \hat{\bm{x}}_k =
    \hat{\bm{\Phi}} \hat{\bm{\Lambda}}^q \hat{\bm{c}}_k,
    \quad
    q=0,\dots,s.
    \label{eq:xkjmatrix}
\end{equation}

\subsection{Augmentation with one sample} \label{subsec:s=1}

When $s=1$, $\hat{\bm{x}}_{k} \in \mathbb{R}^{2n} $ can be written as $ \hat{\bm{x}}_{k} = [\bm{x}_{k}, \bm{x}_{k+1}]^T$, splitting it to the top and bottom $n$ entries. Consequently, \eqref{eq:xkmatrix} can be rewritten as
\begin{equation}
    \begin{bmatrix}
    \bm{x}_{k} \\ \bm{x}_{k+1}
    \end{bmatrix}
    =
    \begin{bmatrix}
    \hat{\bm{\Phi}}^{(0)} \\ \hat{\bm{\Phi}}^{(1)} 
    \end{bmatrix}
    \hat{\bm{c}}_k,
    \label{eq:xkaug}
\end{equation}
where the columns of $\hat{\bm{\Phi}}^{(0)} \in \mathbb{C}^{n \times \hat{r}}$ and $\hat{\bm{\Phi}}^{(1)} \in \mathbb{C}^{n \times \hat{r}}$ are the top and bottom $n$ entries of the columns of $\hat{\bm{\Phi}}$. We term each column of $\hat{\bm{\Phi}}^{(0)}$ or $\hat{\bm{\Phi}}^{(1)}$ as a \emph{DMD sub-mode} or, simply, \emph{sub-mode}. Considering only the top $n$ entries in Eq. \eqref{eq:xkaug} yields

\begin{equation}
    \bm{x}_{k} = \hat{\bm{\Phi}}^{(0)} \hat{\bm{c}}_k.
    \label{eq:xk_phi0_ck}
\end{equation}

On the one hand, Eq. \eqref{eq:xk_phi0_ck} shows that $\bm{x}_{k}$ can be expressed by the $\hat{r}$ columns of $\hat{\bm{\Phi}}^{(0)}$. On the other hand, according to \eqref{eq:xkmatrix_unaugmented}, it can be represented using only the $r \leq \hat{r}$ columns of $\bm{\Phi}$. Therefore, Eq. \eqref{eq:xk_phi0_ck} might not be a compact representation of $\bm{x}_{k}$, leading to the conjecture that the application of delay-coordinates DMD results in two types of DMD components -- those that describe the dynamical system (true) and those that are not related to it and are an artifact of the augmentation (spurious). By applying similar considerations to the bottom $n$ entries in \eqref{eq:xkaug}, the same conjecture can be made for $\bm{x}_{k+1}$ and  $\hat{\bm{\Phi}}^{(1)}$.
In the following assumption, we make these conjectures more precise. We note that this assumption is supported by an empirical verification in Section \ref{sec:results}.

\begin{assumption}
Without loss of generality, we write the matrices $\hat{\bm{\Phi}}^{(j)} \in \mathbb{C}^{n \times \hat{r}}$ as $\hat{\bm{\Phi}}^{(j)} = \left[ \hat{\bm{\Phi}}^{(j)}_{\text{true}}, \hat{\bm{\Phi}}^{(j)}_{\text{spurious}} \right]$, $j=0,1$, where the $r$ leftmost columns $\hat{\bm{\Phi}}^{(j)}_{\text{true}} \in \mathbb{C}^{n \times r}$ are termed \emph{true} and the remaining $\hat{r}-r$ columns $\hat{\bm{\Phi}}^{(j)}_{\text{spurious}} \in \mathbb{C}^{n \times [\hat{r}-r]}$ are termed \emph{spurious}. Assume that the matrix $\hat{\bm{\Phi}}^{(j)}_{\text{true}}$ has full column rank $r$ and that its column space spans the space of observations $\{\bm{x}_k\}$.
    \label{assump:phi0_full_rank}
\end{assumption}
In other words, we assume that the $r$ leftmost columns of $\hat{\bm{\Phi}}^{(0)}$ are the ones required for the compact reduced-ordered representation of $\bm{x}_{k}$ (hence, true), while the rest $\hat{r}-r$ columns are redundant (hence, spurious). 
According to \eqref{eq:xkmatrix_unaugmented} and by Assumption \ref{assump:phi0_full_rank}, implying that the columns of $\hat{\bm{\Phi}}^{(0)}_{\text{true}}$ are linearly independent, $\bm{x}_{k}$ can be represented using only the $r$ true columns of $\hat{\bm{\Phi}}^{(0)}$ and their respective eigenvalues.

Now, consider $\hat{\bm{x}}_{k+1} \in \mathbb{R}^{2n}$, which for $s=1$ is $ \hat{\bm{x}}_{k+1} = [\bm{x}_{k+1}, \bm{x}_{k+2}]^T$. Then, from \eqref{eq:xk1matrix} and by splitting $\hat{\bm{\Phi}}$, $\hat{\bm{\Lambda}}$ and $\hat{\bm{c}}_k$ to their true and spurious parts as above, we have
\begin{equation}
    \begin{bmatrix}
    \bm{x}_{k+1} \\ \bm{x}_{k+2}
    \end{bmatrix}
    =
   \begin{bmatrix}
    \hat{\bm{\Phi}}^{(0)}_{\text{true}}, \hat{\bm{\Phi}}^{(0)}_{\text{spurious}} \\ \hat{\bm{\Phi}}^{(1)}_{\text{true}}, \hat{\bm{\Phi}}^{(1)}_{\text{spurious}}
    \end{bmatrix}
    \begin{bmatrix}
    \hat{\bm{\Lambda}}_{\text{true}} & 0 \\ 0 & \hat{\bm{\Lambda}}_{\text{spurious}}
    \end{bmatrix}
    \begin{bmatrix}
    \hat{\bm{c}}_k^{\text{true}} \\ \hat{\bm{c}}_k^{\text{spurious}}
    \end{bmatrix},
    \label{eq:xk1aug}
\end{equation}
where $\hat{\bm{\Lambda}}_{\text{true}} = \text{diag} [ \hat{\lambda}_1, \dots,  \hat{\lambda}_r]$, $\hat{\bm{\Lambda}}_{\text{spurious}} = \text{diag} [ \hat{\lambda}_{r+1}, \dots,  \hat{\lambda}_{\hat{r}}]$, and $\hat{\bm{c}}_k^{\text{true}}$,  $\hat{\bm{c}}_k^{\text{spurious}}$ denote the $r$ and $\hat{r}-r$ expansion coefficients that correspond to the true and spurious components, respectively.

Under Assumption \ref{assump:phi0_full_rank}, $\bm{x}_{k+1}$ can be represented using only the true parts of \eqref{eq:xkaug} and \eqref{eq:xk1aug} as
\begin{equation}
    \bm{x}_{k+1} = \hat{\bm{\Phi}}^{(1)}_{\text{true}} \hat{\bm{c}}_k^{\text{true}},
    \quad
    \bm{x}_{k+1} = \hat{\bm{\Phi}}^{(0)}_{\text{true}} \hat{\bm{\Lambda}}_{\text{true}} \hat{\bm{c}}_k^{\text{true}}.
    \label{eq:xk1_true}
\end{equation}
Consequently, by equating the right-hand side (RHS) terms of Eq. \eqref{eq:xk1_true}, we have
\begin{equation}
    \left[
    \hat{\bm{\Phi}}_{\text{true}}^{(1)} - \hat{\bm{\Phi}}_{\text{true}}^{(0)} \hat{\bm{\Lambda}}_{\text{true}}
    \right] \hat{\bm{c}}_k^{\text{true}} = \bm{0},
    \quad
    k=0, \dots, m-s-1.
    \label{eq:phic_0}
\end{equation}

\begin{proposition}
For distinct augmented DMD eigenvalues $\hat{\lambda}_l$,
    \begin{equation}
        \hat{\bm{\Phi}}_{\text{true}}^{(1)} = \hat{\bm{\Phi}}_{\text{true}}^{(0)} \hat{\bm{\Lambda}}_{\text{true}}.
        \label{eq:phi1_lambda_phi0_matrix}
    \end{equation}
    \label{prop:phi1_lambda_phi0}
\end{proposition}
Note that assuming distinct DMD eigenvalues is a common practice \cite{bagheri2010analysis,rowley2009spectral,le2017higher}.

\begin{proof}

Denote $\hat{\bm{C}} \in \mathbb{C}^{r \times [m-s]}$ as a matrix whose columns are $\hat{\bm{c}}_k^{\text{true}} $, $k=0, \dots, m-s-1$. Based on the definition of $\hat{\bm{c}}_k$ in Eq. \eqref{eq:xkmatrix}, $\hat{\bm{C}}$ can be recast as a product of a diagonal matrix $\hat{\bm{\eta}}_0 \in \mathbb{C}^{r \times r}$ and a $[m-s]$th order Vandermonde matrix $\hat{\bm{V}}^{m-s}_{\text{and}}(\hat{\lambda}_1, \dots, \hat{\lambda}_r) \in \mathbb{C}^{r \times [m-s]} $, both consisting of the $r$ true DMD eigenvalues, as

\begin{equation}
\begin{split}
    \hat{\bm{C}}
    & =
    \begin{bmatrix}
    | & | &  & | \\
    \hat{\bm{c}}_0^{\text{true}}  & \hat{\bm{c}}_1^{\text{true}} & \cdots & \hat{\bm{c}}_{m-s-1}^{\text{true}} \\
    | & | &  & |
    \end{bmatrix}
    \\
    & =
    \begin{bmatrix}
    \hat{\eta}_{0,1} & 0 & \cdots & 0 \\
    0 & \hat{\eta}_{0,2} & \cdots & 0 \\
    0 & 0 & \ddots & 0 \\
    0 & 0 & \cdots & \hat{\eta}_{0,r}
    \end{bmatrix}
    \begin{bmatrix}
    1 & \hat{\lambda}_1 & (\hat{\lambda}_1)^2 & \cdots & (\hat{\lambda}_1)^{m-s-1} \\
    1 & \hat{\lambda}_2 & (\hat{\lambda}_2)^2 & \cdots & (\hat{\lambda}_2)^{m-s-1} \\
    \vdots & \vdots & \vdots & \ddots & \vdots \\
    1 & \hat{\lambda}_r & (\hat{\lambda}_r)^2 & \cdots & (\hat{\lambda}_r)^{m-s-1}
    \end{bmatrix}
    \\
    & =
    \hat{\bm{\eta}}_0 \hat{\bm{V}}^{m-s}_{\text{and}}(\hat{\lambda}_1, \dots, \hat{\lambda}_r).
\end{split}
\end{equation}
For distinct DMD eigenvalues, the rank of $\hat{\bm{V}}^{m-s}_{\text{and}}(\hat{\lambda}_1, \dots, \hat{\lambda}_r)$ is $\min(r,m-s)$ \cite{pan2020structure}.
Since $r \leq m-s$ (see Eq. \eqref{eq:s_range}), $ \text{rank} \left(  \hat{\bm{V}}^{m-s}_{\text{and}}(\hat{\lambda}_1, \dots, \hat{\lambda}_r) \right) = r $. 
For a non-trivial representation of the dynamical system, i.e., $\hat{\eta}_{0,l} \neq 0$, $l=1, \dots r$, we have that $ \text{rank} \left( \hat{\bm{\eta}}_0 \right) = r $. Consequently, $\text{rank} ( \hat{\bm{C}} )=r$.

For convenience, we define $\bm{\Delta} = \hat{\bm{\Phi}}_{\text{true}}^{(1)} - \hat{\bm{\Phi}}_{\text{true}}^{(0)} \hat{\bm{\Lambda}}_{\text{true}} $.
As the columns of $\hat{\bm{C}}$ are in the subspace of the solutions to Eq. \eqref{eq:phic_0}, then $\dim \left( \ker\left(  \bm{\Delta} \right) \right) = r$.
By the rank-nullity theorem, the rank of the matrix $\bm{\Delta}$ equals to the matrix's number of columns $r$ minus the dimension of its null space, i.e., $\dim \left( \ker\left(  \bm{\Delta} \right) \right)=r$. Therefore, $\text{rank} \left(  \bm{\Delta} \right) = 0$, implying that $\bm{\Delta} = \bm{0}$. Hence, $ \hat{\bm{\Phi}}_{\text{true}}^{(1)} = \hat{\bm{\Phi}}_{\text{true}}^{(0)} \hat{\bm{\Lambda}}_{\text{true}}$, yielding Eq. \eqref{eq:phi1_lambda_phi0_matrix}.
\end{proof}

Denote the columns of $\hat{\bm{\Phi}}_{\text{true}}^{(0)}$ and $\hat{\bm{\Phi}}_{\text{true}}^{(1)}$ by $ \hat{\bm{\phi}}_l^{(0)}$ and $\hat{\bm{\phi}}_l^{(1)}$, respectively.
Their consideration along with $\hat{\lambda}_l$ in $\hat{\bm{\Lambda}}_{\text{true}}$ in \eqref{eq:phi1_lambda_phi0_matrix}, yields 
\begin{equation}
\hat{\bm{\phi}}_l^{(1)} = \hat{\bm{\phi}}_l^{(0)} \hat{\lambda}_l,
\quad
l=1, \dots, r.
\label{eq:phi1_lambda_phi0}
\end{equation}

Proposition \ref{prop:phi1_lambda_phi0} and Eq. \eqref{eq:phi1_lambda_phi0} show the spatiotemporal coupling in delay-coordinates DMD for the special case of $s=1$. In this case, the sub-modes $\hat{\bm{\phi}}_l^{(0)}, \hat{\bm{\phi}}_l^{(1)} \in \mathbb{C}^n$ that comprise the true, augmented mode $\hat{\bm{\phi}}_l \in \mathbb{C}^{2n}$ are related to one another through the corresponding true eigenvalue $\hat{\lambda}_l$ of $\hat{\bm{A}}$. The general relations for augmentation with more than one sample are presented in the next subsection.

\subsection{Augmentation with several samples} \label{subsec:s>1}

In case $s>1$, two observation vectors that are augmented $s$ times, denoted by $\hat{\bm{x}}_{k}, \hat{\bm{x}}_{k+q} \in  \mathbb{R}^{n(s+1)}$, are considered. Then, from Eq. \eqref{eq:xkmatrix} and \eqref{eq:xkjmatrix}, we have
\begin{equation}
\begin{split}
    \hat{\bm{x}}_{k}
    & =
    \begin{bmatrix}
    \bm{x}_{k} \\ \vdots \\ \bm{x}_{k+j} \\ \vdots \\ \bm{x}_{k+s}
    \end{bmatrix}
     =
    \begin{bmatrix}
    \hat{\bm{\Phi}}^{(0)} \\ \vdots \\ \hat{\bm{\Phi}}^{(j)} \\ \vdots \\ \hat{\bm{\Phi}}^{(s)} 
    \end{bmatrix}
    \hat{\bm{c}}_k,
    \\
    \hat{\bm{x}}_{k+q}
    & =
    \begin{bmatrix}
    \bm{x}_{k+q} \\ \vdots \\ \bm{x}_{k+q+i} \\ \vdots \\ \bm{x}_{k+q+s}
    \end{bmatrix}
     =
    \begin{bmatrix}
    \hat{\bm{\Phi}}^{(0)} \\ \vdots \\ 
    \hat{\bm{\Phi}}^{(i)} \\ \vdots \\\hat{\bm{\Phi}}^{(s)} 
    \end{bmatrix}
    \hat{\bm{\Lambda}}^q
    \hat{\bm{c}}_k,
\end{split}
    \label{eq:xk_xkq}
\end{equation}
where the matrix $\hat{\bm{\Phi}}$ in Eq. \eqref{eq:xkmatrix} and \eqref{eq:xkjmatrix} is viewed as a column-stack of $s+1$ ``sub-modes" matrices $\hat{\bm{\Phi}}^{(0)}$, \dots, $\hat{\bm{\Phi}}^{(s)}$, each of dimensions $n \times \hat{r}$.
Without loss of generality, we assume that $j=q+i$, $i,j=0, \dots, s$. Thus, both $\hat{\bm{x}}_{k}$ and $\hat{\bm{x}}_{k+q}$ in \eqref{eq:xk_xkq} contain the same unaugmented observation $\bm{x}_{k+j}=\bm{x}_{k+q+i}$. 
Similarly to the case of $s=1$, by equating the true components of the RHS terms in \eqref{eq:xk_xkq} that correspond to $\bm{x}_{k+j}$ and $\bm{x}_{k+q+i}$ and substituting $j=q+i$, we have
\begin{equation}
    \left[
    \hat{\bm{\Phi}}_{\text{true}}^{(j)} -  \hat{\bm{\Phi}}_{\text{true}}^{(i)} \hat{\bm{\Lambda}}^{j-i}_{\text{true}}
    \right] \hat{\bm{c}}_k^{\text{true}} = \bm{0},
    \quad
    k=0, \dots, m-s-1.
\end{equation}
Next, we generalize Assumption \ref{assump:phi0_full_rank} and Proposition \ref{prop:phi1_lambda_phi0} for augmentation with more than one sample, where the spatiotemporal coupling in delay-coordinates DMD is bounded from above by $B$, which is defined and explained in the sequel in Eq. \eqref{eq:bound}.

\begin{assumption}
    The matrices $\hat{\bm{\Phi}}^{(j)}_{\text{true}}$, $j=0, \dots, B$, comprising the $r$ leftmost columns of $\hat{\bm{\Phi}}^{(j)}$, are assumed to have full column rank and that the column space of each of them spans the space of observations $\{\bm{x}_k\}$.
    \label{assump:phij_full_rank}
\end{assumption}

\begin{figure*}[t]
\includegraphics[width=0.7\textwidth]{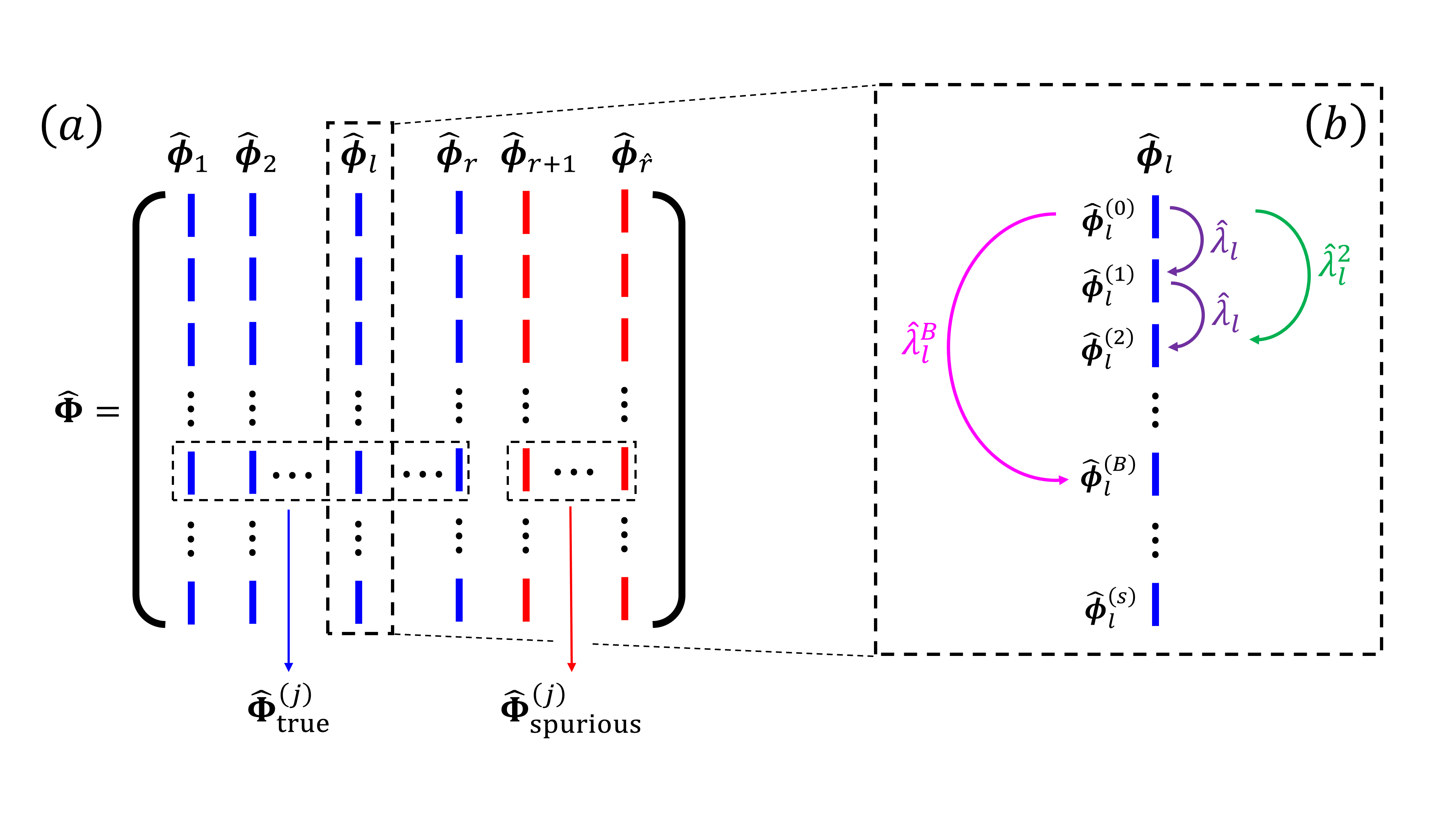}
\caption{
\label{fig:sdcdmd_illustration}
(a) Division of $\hat{\bm{\Phi}}$, whose columns are the augmented DMD modes $\hat{\bm{\phi}}_l$, into $\hat{\bm{\Phi}}_{\text{true}} \in \mathbb{C}^{[n(s+1)] \times r}$ and $\hat{\bm{\Phi}}_{\text{spurious}} \in \mathbb{C}^{[n(s+1)] \times [\hat{r}-r]}$. The columns of the former and latter are the $r$ true and $\hat{r}-r$ spurious augmented DMD modes, respectively.
(b) Illustration of the spatiotemporal coupling in delay-coordinates DMD, according to Proposition \ref{prop:phistruct} and Eq. \eqref{eq:phistruct}.
}
\end{figure*}

The division of $\hat{\bm{\Phi}}$ into $\hat{\bm{\Phi}}^{(j)}_{\text{true}} \in \mathbb{C}^{n \times r}$ and $\hat{\bm{\Phi}}^{(j)}_{\text{spurious}} \in \mathbb{C}^{n \times [\hat{r}-r]}$, $j=0, \dots, B$, as in Assumptions \ref{assump:phi0_full_rank} and \ref{assump:phij_full_rank}, is illustrated in Fig. \ref{fig:sdcdmd_illustration}(a).

\begin{proposition}
    For distinct augmented DMD eigenvalues $\hat{\lambda}_l$,
    \begin{equation}
    \hat{\bm{\Phi}}_{\text{true}}^{(j)} = \hat{\bm{\Phi}}_{\text{true}}^{(i)} \hat{\bm{\Lambda}}^{j-i}_{\text{true}},
    \quad
    i,j=0, \dots, B.
    \label{eq:phistruct_matrix}
    \end{equation}
    \label{prop:phistruct}
\end{proposition}
We omit the proof because it is similar to the proof of Proposition \ref{prop:phi1_lambda_phi0}.

Consideration of $\hat{\bm{\phi}}_l^{(i)}$ and $\hat{\bm{\phi}}_l^{(j)}$, which are the respective $i$th and $j$th sub-modes of the true, augmented modes $\hat{\bm{\phi}}_l$  (i.e., the columns of  $\hat{\bm{\Phi}}_{\text{true}}^{(i)}$ and $\hat{\bm{\Phi}}_{\text{true}}^{(j)}$ in \eqref{eq:phistruct_matrix}, respectively), as well as $\hat{\lambda}_l$ in $\hat{\bm{\Lambda}}_{\text{true}}$, $l=1, \dots, r$, leads to
\begin{equation}
    \hat{\bm{\phi}}_l^{(j)} = \hat{\bm{\phi}}_l^{(i)} \hat{\lambda}_l^{j-i},
    \quad
    i,j=0, \dots, B,
    \quad
    l=1, \dots, r.
    \label{eq:phistruct}
    \end{equation}
    
The spatiotemporal coupling in delay-coordinates DMD, generally presented in Proposition \ref{prop:phistruct} and Eq. \eqref{eq:phistruct}, reveals that the true, augmented modes $\hat{\bm{\phi}}_l$ contain not only spatial information, but also temporal information given by $\hat{\lambda}_l$.
Importantly, by Proposition \ref{prop:phistruct}, this coupling is exhibited only by the true DMD components. In contrast, the considerations leading to this property do not apply to the spurious DMD components, and therefore, they lack this coupling.
This spatiotemporal coupling can be written more explicitly using \eqref{eq:phistruct}, e.g., as

\begin{equation}
\begin{split}
    \hat{\bm{\phi}}_l &=
    \begin{bmatrix}
    \hat{\bm{\phi}}_l^{(0)} \\ \hat{\bm{\phi}}_l^{(1)} \\ \hat{\bm{\phi}}_l^{(2)} \\ \vdots \\ \hat{\bm{\phi}}_l^{(j)}
    \end{bmatrix}
    =
    \begin{bmatrix}
    \hat{\bm{\phi}}_l^{(0)} \\ \hat{\bm{\phi}}_l^{(0)} \hat{\lambda}_l \\  \hat{\bm{\phi}}_l^{(0)} \hat{\lambda}_l^2 \\ \vdots \\ 
    \hat{\bm{\phi}}_l^{(0)} \hat{\lambda}_l^j
    \end{bmatrix},
    \\
    j &= 0, \dots, B,
    \\
    l&=1, \dots, r,
\end{split}
    \label{eq:explicitstructure}
\end{equation}
and is illustrated in Fig. \ref{fig:sdcdmd_illustration}(b), showing a zoom-in on the single augmented DMD mode $\hat{\bm{\phi}}_l$.
Alternatively, one can obtain the augmented DMD eigenvalue through \eqref{eq:phistruct} by
\begin{equation}
\begin{split}
    \hat{\lambda}_l & = \left( \frac{\hat{\bm{\phi}}_l^{(j)}}{\hat{\bm{\phi}}_l^{(i)}} \right) ^ {\frac{1}{{j-i}}},
    \\
    i,j &= 0, \dots, B, \ i \neq j,
    \\
    l &= 1,\dots,r.
\end{split}
    \label{eq:phistrcturegeneral}
\end{equation}

Thus far, we claimed that the spatiotemporal coupling in Eqs. \eqref{eq:phistruct_matrix} - \eqref{eq:phistrcturegeneral} holds for augmentation numbers bounded from above by $B$. Now, we present this upper bound explicitly and show that it depends on the oscillation frequencies of the underlying continuous-time dynamical system, $\omega_l$. The existence of the upper bound $B$ stems from the Nyquist-Shannon sampling criterion\cite{nyquist1928certain,shannon1949communication}, provided that the sampling frequency $\omega_s$ is sufficiently high, i.e., $\omega_l \leq 0.5 \omega_s$ for every $\omega_l$ of the system.

By substituting $\Delta t = 2 \pi / \omega_s$ into \eqref{eq:omegafromlambda}, the relationship between the DMD eigenvalue $\hat{\lambda}_l$ and the oscillation frequency $\omega_l$ is
\begin{equation}
    \hat{\lambda}_l = \exp(\omega_l \Delta t) = \exp \left( \frac{\omega_l}{0.5 \omega_s} \pi \right).
\end{equation}
Therefore, taking powers $p \in \mathbb{R}$ of $\hat{\lambda}_l$ results in
\begin{equation}
    \hat{\lambda}_l^p
    =
    \left( \exp \left( \frac{\omega_l}{0.5 \omega_s} \pi \right) \right)^p
    =
    \exp \left(\frac{p \omega_l}{0.5 \omega_s} \pi \right).
    \label{eq:lambdap}
\end{equation}
Namely, $\hat{\lambda}_l^p$ corresponds, in effect, to a frequency that is $p$-times faster than the continuous time counterpart of $\hat{\lambda}_l$, i.e., $p \omega_l$. So, the Nyquist-Shannon sampling criterion for this frequency is $p \omega_l \leq 0.5 \omega_s$, which provides the constraint $p \le 0.5\omega_s/\omega_l$. Accordingly, we define
\begin{equation}
   B(\omega_l) = \min\left(s, \left\lfloor \frac{0.5 \omega_s}{\omega_l} \right\rfloor \right),
    \label{eq:bound}
\end{equation}
where $\lfloor \cdot \rfloor$ is the floor function.

In terms of Eq. \eqref{eq:explicitstructure} and as illustrated in Fig. \ref{fig:sdcdmd_illustration} (b), $\hat{\lambda}_l^j$ relates the zeroth sub-mode of $\hat{\bm{\phi}}_l$ to its $j$th sub-mode, where $\hat{\bm{\phi}}_l$ represents the $l$th DMD mode that corresponds to the oscillation frequency $\omega_l$.
Additionally, according to Eq. \eqref{eq:lambdap}, $\hat{\lambda}_l^j$ represents a $j$-times faster frequency than $\hat{\lambda}_l$, and is utilized by the spatiotemporal coupling in delay-coordinates DMD in Proposition \ref{prop:phistruct}.
Consequently, this spatiotemporal coupling exists only for values of $j$ that satisfy the Nyquist-Shannon sampling criterion (i.e., $j=0, \dots, B$). In other words, $\hat{\lambda}_l^j$ for $j>0.5 \omega_s / \omega_l$ represents a frequency that is sampled at a sub-Nyquist rate, and therefore, does not capture the underlying dynamics.

We showed that the sub-modes within an augmented mode are related to each other through the corresponding eigenvalue (see Proposition \ref{prop:phistruct}, Eq. \eqref{eq:phistruct} and Fig. \ref{fig:sdcdmd_illustration}(b)). We refer to these relations as the \textit{spatiotemporal coupling in delay-coordinates DMD}. That is, an augmented DMD mode holds temporal information in addition to the spatial information. In addition, we show empirically in Section \ref{sec:example} that the spurious components do not have this coupling.
In the following sections we utilize the spatiotemporal coupling for selection of the true DMD components and further representation and characterization of dynamical systems from observations, even when the observations are corrupted with noise.


\section{Proposed method}\label{sec:proposed_method}

Based on the spatiotemporal coupling in delay-coordinates DMD presented in Proposition \ref{prop:phistruct}, we propose a method for obtaining a compact and reduced-order representation of the observations that is analogous to Eq. \eqref{eq:xk}. 
As discussed in Section \ref{sec:sdcdmd}, representation \eqref{eq:xk} cannot be directly obtained from applications of DMD algorithms to delay-coordinates embedding due to two challenges: (a) the possible existence of both true and spurious augmented DMD components, and (b) the existence of redundant sub-modes comprising every augmented DMD mode, which are unnecessary for the representation. Both challenges must be addressed to obtain the desired compact representation, i.e., the selection of the true augmented DMD modes, followed by the selection of the relevant sub-modes within the true modes.

To the best of our knowledge, challenge (a) is not addressed in the current literature.
A common practice to address challenge (b) is to choose the first $n$ entries in each selected augmented DMD mode, ignoring its spatiotemporal coupling.
Therefore, we propose a method that considers the spatiotemporal coupling in Proposition \ref{prop:phistruct} and selects the DMD components required for compact representations of dynamical systems. Our method is detailed in this section and summarized in Algorithm \ref{alg}. An empirical comparison between our method and the maximal amplitudes method is conducted in Section \ref{sec:results}, where the superiority of our method is demonstrated in case the observations are corrupted with high levels of noise.

Our method begins by augmenting the observations $s$ times, and applying the exact DMD to them. Generally, the augmentation number $s$ can be chosen either in the range prescribed in \eqref{eq:s_range_small} or in \eqref{eq:s_range}. Yet, the empirical evidence presented in Section \ref{sec:results} shows that $s$ should be chosen according to \eqref{eq:s_range}.
Next, the DMD eigenvalue-mode pairs, $\hat{\lambda}_l$ and $\hat{\bm{\phi}}_l$, respectively, are extracted, followed by a calculation of the oscillation frequencies $\omega_l$ by \eqref{eq:omegafromlambda}, as well as the bounds $B(\omega_l)$ that correspond to each $\omega_l$ by \eqref{eq:bound}.

To identify the true components out of all the obtained DMD components, we propose to utilize Eq. \eqref{eq:phistrcturegeneral}, where the identification proceeds as follows.
The same eigenvalue $\Tilde{\lambda}_l$ is repeatedly computed via substitution of different $i$ and $j$ values in Eq. \eqref{eq:phistrcturegeneral}. For example, by setting $i=0$ and $j=1, \dots , B(\omega_l)$, the eigenvalue $\Tilde{\lambda}_l$ can be computed multiple times by taking the quotients of the sub-modes as
\begin{equation}
    \Tilde{\lambda}_l^{(j)} = \left( \frac{\hat{\phi}_l^{(j)}}{\hat{\phi}_l^{(0)}} \right) ^ {\frac{1}{j}}, \quad j=1, \dots, B(\omega_l),
    \label{eq:computed_lambda}
\end{equation}
where a tilde denotes a computed entity, and the superscripted index $(j)$ in $\Tilde{\lambda}_l^{(j)}$ denotes the $j$th computation of this eigenvalue.
As stated in Proposition \ref{prop:phistruct}, the spatiotemporal coupling in delay-coordinates DMD holds only for true DMD components. Hence, any $j$th computation of any $l$th true eigenvalue according to \eqref{eq:computed_lambda} yields the same value (in a noise-free system).
Moreover, this value is identical to the DMD eigenvalue, i.e., $\Tilde{\lambda}_l^{(j)} = \hat{\lambda}_l$, for $l=1,\dots,r$.
Conversely, such computations of spurious DMD eigenvalues, for $l=r+1,\dots,\hat{r}$, are excepted to yield different computed results upon substitution of varying values of $j$ in \eqref{eq:computed_lambda}, as well as results that differ from $\hat{\lambda}_l$.
Therefore, \eqref{eq:computed_lambda} contains information that enables the distinction between true and spurious DMD components.

In the presence of observation noise, the identification of true DMD components through \eqref{eq:computed_lambda} can be enhanced by introducing an averaged computed eigenvalue, $\langle \Tilde{\lambda}_l \rangle$, e.g., as
\begin{equation}
    \langle \Tilde{\lambda}_l \rangle = \frac{1}{B(\omega_l)} \sum_{j=1}^{ B(\omega_l) } \Tilde{\lambda}_l^{(j)}.
    \label{eq:avg_lambda}
\end{equation}
In case $\omega_l>0.5\omega_s$, then $B(\omega_l)=0$.
This scenario can occur, e.g., when $s$ is chosen to be very large, which, in turn, might produce a very large $\omega_l$. Since the sampling rate admits the Nyquist-Shannon sampling criterion for the true oscillation frequencies, such $\omega_l$ must be related to a spurious eigenvalue. Therefore, when $B(\omega_l)=0$, we consider the eigenvalue corresponding to $\omega_l$ as spurious and remove it from the computation in \eqref{eq:avg_lambda}.

Thereafter, the absolute errors, $\varepsilon_l$, are calculated by
\begin{equation}
    \varepsilon_l = \left| \langle \Tilde{\lambda}_l \rangle - \hat{\lambda}_l \right|,
    \label{eq:eig_error}
\end{equation}
providing an estimation for the difference between the augmented DMD eigenvalue and the average of the eigenvalues computed based on the spatiotemporal coupling in Proposition \ref{prop:phistruct}.
Alternatively, the absolute errors can be defined as the average of the absolute differences between computed eigenvalues $\Tilde{\lambda}_l^{(j)}$ and their corresponding DMD eigenvalue $\hat{\lambda}_l$. We note that this definition yields similar empirical results in the applications studied in this paper.

If the computed eigenvalues $\Tilde{\lambda}_l^{(j)}$ for $j=1, \dots, B(\omega_l)$ are different from $\hat{\lambda}_l$, then $ \varepsilon_l$ in Eq. \eqref{eq:eig_error} is large, which means that $\hat{\lambda}_l$ is spurious.
Contrarily, if $\Tilde{\lambda}_l^{(j)}$ are approximately equal to $\hat{\lambda}_l$, then $\varepsilon_l$ is small, indicating that $\hat{\lambda}_l$ is true.
Therefore, $\varepsilon_l$, $l=1, \dots, \hat{r}$, can be partitioned into two subsets, $S_1$ and $S_2$. Without loss of generality, assume that $S_1$ constitutes the first $d$ values of $\varepsilon_l$, i.e., $\{ \varepsilon_1, \dots, \varepsilon_d \} \in S_1$, where $1 \leq d \leq \hat{r}$, and that $S_2$ contains the remaining $\hat{r}-d$ values, $ \{ \varepsilon_{d+1}, \dots, \varepsilon_{\hat{r}} \} \in S_2$. This partitioning can be carried out by applying standard clustering algorithms to $\varepsilon_l$, e.g., $k$-means with $k=2$.
After clustering, the average absolute errors of subsets $S_1$ and $S_2$, $ \langle \varepsilon \rangle _{S_1}$ and $ \langle \varepsilon \rangle _{S_2}$ are computed as
\begin{equation}
    \langle \varepsilon \rangle _{S_1} = \frac{1}{d}\sum_{l=1}^{d} \varepsilon_l,
    \quad
    \langle \varepsilon \rangle _{S_2} = \frac{1}{\hat{r}-d}\sum_{l=1}^{\hat{r}-d} \varepsilon_l,
    \label{eq:avg_eig_error}
\end{equation}
and then, the smaller of $\langle \varepsilon \rangle _{S_1}$ and $\langle \varepsilon \rangle _{S_2}$ denotes the subset that contains the $\varepsilon_l$ related to the true DMD eigenvalues.
The distinction between small and large absolute errors was straightforward in the illustrative example (see Section \ref{sec:example}, Fig. \ref{fig:abs_err_example}).

In some cases, the DMD components of dynamical systems arise in complex conjugate pairs.
In such cases, each complex conjugate pair represents the same characteristics of the system (e.g., the same oscillation frequencies $\omega_l$). These pairs embody another degree of redundancy, which can be exploited, for example, as follows.
Denote $\hat{\lambda}_q$ as the complex conjugate of $\hat{\lambda}_l$. Then, similarly to \eqref{eq:avg_lambda}, an averaged eigenvalue $\langle \Tilde{\lambda}_{l,q} \rangle$, which accounts for this redundancy, can be written as
\begin{equation}
    \langle \Tilde{\lambda}_{l,q} \rangle = \frac{1}{2} \left( \langle \Tilde{\lambda}_l \rangle  + \overline{ \langle \Tilde{\lambda}_q \rangle } \right),
    \label{eq:avg_lambda_complex_conjugate}
\end{equation}
where an overline denotes a complex conjugate.
As $\langle \Tilde{\lambda}_{l,q} \rangle$ is computed over both $l$ and $q$, the influence of noise on the observations can be further diminished in these systems.

Lastly, once the subset of the true augmented DMD components is identified using Eq. \eqref{eq:avg_eig_error}, the observations can be represented analogously to \eqref{eq:xk} as
\begin{equation}
    \bm{x}_k = \sum_{l=1}^{|S_1|}  \langle \Tilde{\lambda}_l \rangle ^k \hat{\bm{\phi}}_l^{(0)} \hat{\sigma}_{0,l},
    \label{eq:xk_analogous}
\end{equation}
where $\hat{\sigma}_{0,l}=\langle \hat{\bm{\phi}}_l,\bm{\hat{x}}_0 \rangle $, $|S_1|$ is the cardinality of $S_1$, and, without loss of generality, we assume that $S_1$ is the subset that contains the $\varepsilon_l$ values related to the true DMD eigenvalues $\langle \Tilde{\lambda}_l \rangle$.

We note that our method has two main shortcomings, which we plan to address in future research.
First, compared to mode selection based on maximal amplitudes, our method is computationally heavier. Second, as shown in the next section, our empirical study suggests that the performance is sensitive to the choice of the augmentation number $s$. We plan to develop a systematic procedure to set the augmentation number, as well as more efficient implementation schemes that mitigate the repeated computation of eigenvalue decomposition. 


\begin{algorithm}[H]

\caption{Proposed algorithm}\label{alg}
\flushleft \textbf{Input}: Observations $\{\bm{x}_k\}_{k=0}^{m}$ of a dynamical system. \\
\textbf{Output}: Compact and reduced-order spectral representation of the observations.

\begin{enumerate}
    \item Augment the data $s$ times according to the range in \eqref{eq:s_range}, apply the exact DMD algorithm to the augmented observations, and extract the augmented DMD eigenvalue-mode pairs, $\hat{\lambda}_l$ and $\hat{\bm{\phi}}_l$, respectively.

    \item Calculate the oscillations frequencies $\omega_l$ that correspond to each augmented DMD eigenvalue $\hat{\lambda}_l$ using Eq. \eqref{eq:omegafromlambda}.
    
    \item Compute the bound $B(\omega_l)$ in Eq. \eqref{eq:bound} for each $\omega_l$.
    
    \item  Compute the averaged eigenvalues $\langle 
    \Tilde{\lambda}_l \rangle$ by Eq. \eqref{eq:avg_lambda}.
    
    \item Calculate the absolute errors, $\varepsilon_l$, between the averaged computed eigenvalues, $\langle \Tilde{\lambda_l} \rangle$, and their corresponding (non-averaged) DMD eigenvalues, $\hat{\lambda}_l$, as in Eq. \eqref{eq:eig_error}. 
    
    \item Partition $\varepsilon_l$ into two subsets using $k$-means with $k=2$, and compute the average values of each subset, as in \eqref{eq:avg_eig_error}. Select the averaged eigenvalues (and their corresponding modes) related to the subset corresponding to the smaller average value.

    \item Use the selected DMD components for building a compact and reduced-order spectral representation of the observations, according to Eq. \eqref{eq:xk_analogous}.
    
\end{enumerate}
\end{algorithm}


\section{Illustrative example}\label{sec:example}

\begin{figure}[t]
\includegraphics[width=\linewidth]{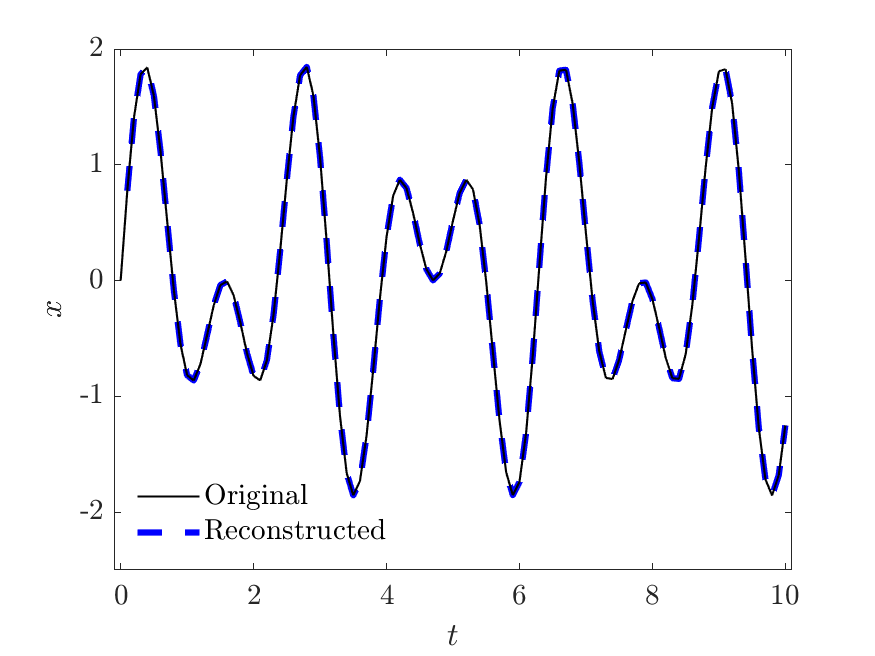}
\caption{
\label{fig:reconstructed}
The two-mode sine signal in Eq. \eqref{eq:twomode}, depicted as the black solid line (original signal). The blue dashed line depicts the reconstructed signal according to Algorithm \ref{alg} and Eq. \eqref{eq:xk_analogous}.
}
\end{figure}

By following the proposed algorithm (Algorithm \ref{alg}) step by step, we demonstrate the spatiotemporal coupling in delay-coordinates DMD (formulated in section \ref{sec:sdcdmd}) on an illustrative example of a two-mode sine signal.
More specifically, we show how Proposition \ref{prop:phistruct} facilitates the identification, in an equation-free manner, of the DMD components required for the representation and characterization of the signal. Evaluation is done by comparing the true signal's oscillation frequencies to the oscillation frequencies that correspond to the DMD eigenvalues selected by our method.
Our method can also be visually evaluated through its signal reconstruction compared to the original signal, as shown in Fig. \ref{fig:reconstructed}.


\begin{figure}

\begin{center}

\begin{tikzpicture}[black!75,thick]
 
\fill [pattern = north west lines] (0,0) rectangle (7,0.5);
\draw[thick] (0,0.5) -- (7,0.5);

\fill [pattern = north west lines] (-0.5,0) rectangle (0,2);
\draw[thick] (0,0.5) -- (0,2);

\draw
[
    decoration={
        coil,
        aspect=0.3, 
        segment length=1.2mm, 
        amplitude=2mm, 
        pre length=3mm,
        post length=3mm},
    decorate
] (0,1.25) -- (2,1.25) 
    node[midway,above=0.25cm,black]{};

\node[draw,
    fill=white!60,
    minimum width=1cm,
    minimum height=1cm,
    anchor=north,
    label=center:$m_1$] at (2.5,1.75) {} ;
    \draw [thick, fill=black] (2.25,0.62) circle (0.1) (2.75,0.62) circle (0.1);

\draw[very thick,
    black,
    |-latex] (2.5,2) -- (3.5,2)
    node[midway,above]{$x(t)$};

\draw
[
    decoration={
        coil,
        aspect=0.3, 
        segment length=1.2mm, 
        amplitude=2mm, 
        pre length=3mm,
        post length=3mm},
    decorate
] (3,1.25) -- (5,1.25) 
    node[midway,above=0.25cm,black]{}; 

\node[draw,
    fill=white!60,
    minimum width=1cm,
    minimum height=1cm,
    anchor=north,
    label=center:$m_2$] at (5.5,1.75) {};
    \draw [thick, fill=black] (5.25,0.62) circle (0.1) (5.75,0.62) circle (0.1);


\end{tikzpicture}
\end{center}

\caption{\label{fig:2dof_undamped_oscillator}
Two degrees of freedom oscillator, comprised of two masses, $m_1$ and $m_2$, connected via two springs.
}

\end{figure}
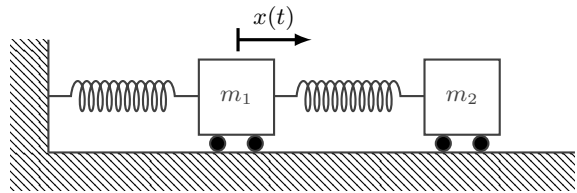


Consider the following two-mode sine signal
\begin{equation}
    x(t)=\sin(\omega_1^\text{sys} t)+\sin(\omega_2^\text{sys} t) + \alpha n(t),
    \label{eq:twomode}
\end{equation}
where $\omega_1^{\text{sys}}=3$ rad/s,  $\omega_2^{\text{sys}}=5$ rad/s, $n(t)$ is an additive white standard Gaussian noise, and $\alpha=10^{-12}$.
Fig. \ref{fig:2dof_undamped_oscillator} shows a dynamical system whose coordinate $x(t)$ can represent such a signal. Specifically, signal \eqref{eq:twomode} can be viewed as the position of mass $m_1$ in a two degrees of freedom (DOF) oscillator, comprised of two masses, $m_1$ and $m_2$, connected to each other via springs.
We denote the components that correspond to the system by the superscript sys, i.e., $\omega_1^{\text{sys}}$ and $\omega_2^{\text{sys}}$ are the system's harmonics.

The signal is sampled at the sampling rate $\omega_s= 2\pi / \Delta t = 20\pi$ rad/s (corresponding to 10 Hz) and 101 samples (observations) $\{x_k\}_{k=0}^{100}$ are collected, shown by the solid black line in Fig. \ref{fig:reconstructed}.

We begin the illustration by augmenting the observations with $s=11$. This choice of $s$ is in accordance with the upper bound on $s$ in \eqref{eq:s_range}, yet it does not adhere to the lower bound in the same equation. One of the main purposes of this section is to demonstrate that Proposition \ref{prop:phistruct} holds for true DMD components and up to the bound $B(\omega_l)$.
Consequently, we choose $s=11$, which is a value that, on the one hand, demonstrates these two statements for system \eqref{eq:twomode}, and on the other hand, provides a small enough number of eigenvalues that can be conveniently visualized. That is, this choice is made merely for illustrative purposes.

\begin{figure}[t]
\includegraphics[width=\linewidth]{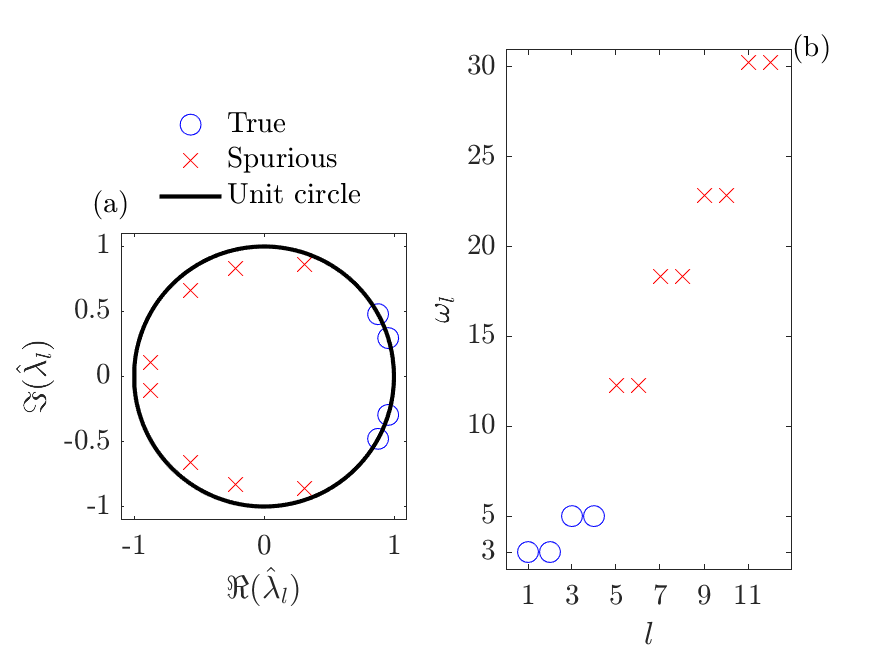}
\caption{
\label{fig:true_spurious}
(a) Polar representation of the augmented DMD eigenvalues, $\hat{\lambda}_l$, that arise from an application of the exact DMD algorithm to the observations $\{\bm{x}_k\}_{k=0}^{100}$ sampled from signal \eqref{eq:twomode}, which were augmented $s=11$ times. The black line denotes the unit circle. (b) The oscillation frequencies, $\omega_l$, that correspond to $\hat{\lambda}_l$ through Eq. \eqref{eq:omegafromlambda}. The \emph{true} eigenvalues that are related to the system's harmonics ($\omega_1^{\text{sys}}=3$ rad/s and $\omega_2^{\text{sys}}=5$ rad/s) are marked by blue circles, while the red crosses mark the \emph{spurious} eigenvalues that are not related to them.
}
\end{figure}

The exact DMD algorithm is applied to the augmented observations, resulting in $s+1=12$ augmented DMD eigenvalue-mode pairs, $\{ \hat{\lambda}_l, \hat{\bm{\phi}}_l \}_{l=1}^{12}$.
Fig. \ref{fig:true_spurious}(a) presents the polar representation of $\hat{\lambda}_l$, and Fig. \ref{fig:true_spurious}(b) shows the oscillation frequencies, $\omega_l$, (corresponding to the eigenvalues by Eq. \eqref{eq:omegafromlambda}), sorted in ascending order.
It can be observed in Fig. \ref{fig:true_spurious}(a) that the DMD eigenvalues of this signal arise in complex conjugate pairs. Moreover, these eigenvalues are distinct, satisfying the condition of Propositions \ref{prop:phi1_lambda_phi0} and \ref{prop:phistruct}.

Since the signal is composed of two oscillation frequencies, $\omega_1^{\text{sys}}$ and $\omega_2^{\text{sys}}$, two complex conjugate DMD eigenvalue pairs $\hat{\lambda}_{1,2}$ and $\hat{\lambda}_{3,4}$ are related to these frequencies $-$ these are the \emph{true} eigenvalues and they are marked by blue circles in Fig. \ref{fig:true_spurious}. Indeed, Fig. \ref{fig:true_spurious}(b) corroborates them as true since they correspond to the frequencies 3 rad/s and 5 rad/s. Similarly, the remaining 8 eigenvalues (and their corresponding modes) are \emph{spurious} and are marked by red crosses in Fig. \ref{fig:true_spurious}.

Our goal of identifying the true DMD components required for the system description in an equation-free manner can be visually described via Fig. \ref{fig:true_spurious}. Suppose the observations $\{x_k\}_{k=0}^{100}$ were obtained from an unknown system, i.e., without the knowledge that 3 rad/s and 5 rad/s are the system's harmonics. Then, the markings of the blue circles and red crosses would be unknown, as well. In this regard, our goal would be the task of identifying the blue circles, which represent the true DMD eigenvalues.

\begin{figure*}[t]
\includegraphics[width=0.7\textwidth]{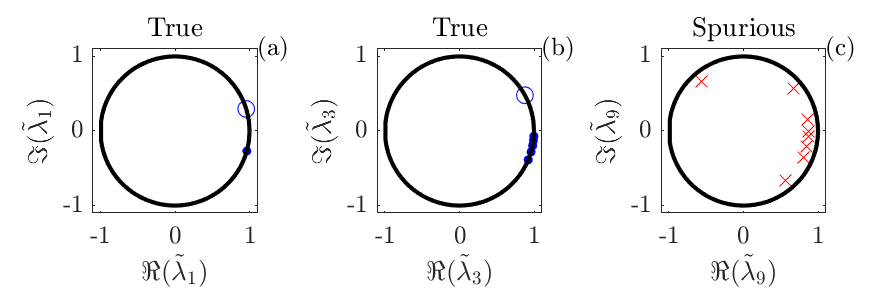}
\caption{
\label{fig:structure_limit}
Polar representation of the computed eigenvalues $\Tilde{\lambda}_1^{(j)}$, $\Tilde{\lambda}_3^{(j)}$, $\Tilde{\lambda}_9^{(j)}$, $j=1, \dots, 11$, obtained from Eq. \eqref{eq:computed_lambda}. The black lines denote the unit circle.
According to the bound $B(\omega_l)$ in Eq. \eqref{eq:bound}, which is valid only for \emph{true} DMD components, the first 10 (a) and 6 (b) values of the \emph{true} $\Tilde{\lambda}_1^{(j)}$ and $\Tilde{\lambda}_3^{(j)}$, respectively, are the same. These values are denoted by blue circles, which coincide with each other 10 times in (a) and 6 times in (b).
Contrarily, the remaining eleventh value $\Tilde{\lambda}_1^{(11)}$ in (a), and the five remaining values $\Tilde{\lambda}_3^{(j)}$, $j=7, \dots, 11$ in (b), 
 which are marked by blue dots, differ from one another and from the true $\hat{\lambda}_1$ or $\hat{\lambda}_3$.
(c) Repeated computations of $\Tilde{\lambda}_9^{(j)}$ that yield different values, since the relations in Proposition \ref{prop:phistruct} do not hold for \emph{spurious} DMD components.
}

\end{figure*}

By Eq. \eqref{eq:bound}, the bounds that correspond to the two system's harmonics are computed to be $B(\omega_1^{\text{sys}}) = 10$ and $B(\omega_2^{\text{sys}}) = 6$. Next, we demonstrate Proposition \ref{prop:phistruct} by showing that the spatiotemporal coupling indeed holds up to these bounds. The choice of $s=11$ is the minimal value of $s$ that can serve this purpose.
As an example, we consider the true eigenvalues $\hat{\lambda}_1$ and $\hat{\lambda}_3$, as well as the spurious eigenvalue $\hat{\lambda}_9$. For each of them, we follow the eigenvalues computation in Eq. \eqref{eq:computed_lambda} eleven times; namely, we substitute $j=1, \dots ,11$ for each of the three values $l=1, 3, 9$ in Eq. \eqref{eq:computed_lambda}. Fig. \ref{fig:structure_limit}
presents the polar representations of these computed eigenvalues.

Fig. \ref{fig:structure_limit}(a) shows $\Tilde{\lambda}_1^{(1)}, \dots, \Tilde{\lambda}_1^{(10)}$ as blue circles and $\Tilde{\lambda}_1^{(11)}$ as a blue dot. As expected and according to $B(\omega_1^{\text{sys}}) = 10$, the former 10 computed eigenvalues coincide with each other, and appear as a single circle in Fig. \ref{fig:structure_limit}(a). Moreover, a comparison between them and their corresponding true DMD eigenvalue, $\hat{\lambda}_1$, reveals that $\Tilde{\lambda}_1^{(1)} = \dots =  \Tilde{\lambda}_1^{(10)} = \hat{\lambda}_1$. Contrarily, the 11\textsuperscript{th} computed eigenvalue, $\Tilde{\lambda}_1^{(11)}$, differs from them, as expected, since its index exceeds the upper bound $B(\omega_1^{\text{sys}}) = 10$.

Similarly, Fig. \ref{fig:structure_limit}(b) shows $\Tilde{\lambda}_3^{(1)}, \dots, \Tilde{\lambda}_3^{(6)}$ as blue circles and $\Tilde{\lambda}_3^{(7)}, \dots, \Tilde{\lambda}_3^{(11)}$ as blue dots. As $B(\omega_2^{\text{sys}}) = 6$, the former 6 computed eigenvalues coincide with each other and appear as a single circle in Fig. \ref{fig:structure_limit}(b). Further comparison confirms that $\Tilde{\lambda}_3^{(1)} = \dots = \Tilde{\lambda}_3^{(6)} = \hat{\lambda}_3$. On the other hand, $\Tilde{\lambda}_3^{(7)} \neq \dots \neq \Tilde{\lambda}_3^{(11)} \neq \hat{\lambda}_3$. That is, relations \eqref{eq:computed_lambda} are indeed valid for true DMD components and up to $B(\omega_l)$.

Lastly, Fig. \ref{fig:structure_limit}(c), which relates to the spurious DMD eigenvalue $\hat{\lambda}_9$, depicts the 11 computed eigenvalues $\Tilde{\lambda}_9^{(1)}, \dots, \Tilde{\lambda}_9^{(11)}$, which are different from each other. This difference indicates that $\hat{\lambda}_9$ is spurious, as relations \eqref{eq:computed_lambda} are not valid for spurious DMD components.

To conclude the demonstration in Fig. \ref{fig:structure_limit}, we showed that only \emph{true} DMD components adhere to the spatiotemporal coupling in delay-coordinates DMD in Proposition \ref{prop:phistruct}, and only up the bound $B(\omega_l)$ in Eq. \eqref{eq:bound}.

\begin{figure}[t]
\includegraphics[width=\linewidth]{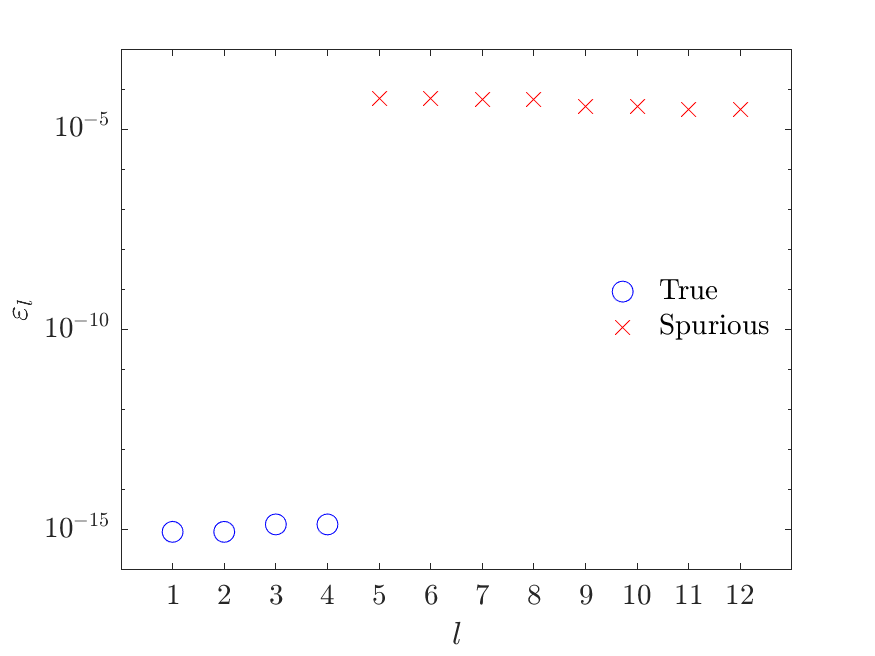}
\caption{
\label{fig:abs_err_example}
The absolute errors, $\varepsilon_l$, between the averaged eigenvalues $\langle \Tilde{\lambda}_l \rangle$ in Eq. \eqref{eq:avg_lambda} and their corresponding DMD eigenvalues, $\hat{\lambda}_l$ (see Eq. \eqref{eq:eig_error}).
The true eigenvalues, $\hat{\lambda}_1, \dots, \hat{\lambda}_4$ are marked by blue circles, and the spurious eigenvalues, $\hat{\lambda}_5, \dots, \hat{\lambda}_{12}$, are marked by red crosses.
The absolute errors of the spurious $\hat{\lambda}_l$ are larger by at least 10 orders of magnitude than those of the true $\hat{\lambda}_l$, providing a clear distinction between the true and spurious DMD eigenvalues.
The markings of the blue circles and red crosses were obtained in a unsupervised manner using Algorithm \ref{alg}.
}
\end{figure}

Next, we exploit this spatiotemporal coupling in order to identify the true DMD components required for compact representation of the dynamical system. Accordingly, we compute all the averaged eigenvalues, $\langle 
\Tilde{\lambda}_l \rangle$, $l=1, \dots, 12$, via Eq. \eqref{eq:avg_lambda}, followed by a calculation of their absolute errors, $\varepsilon_l$, by Eq. \eqref{eq:eig_error}.
Fig. \ref{fig:abs_err_example} shows $\varepsilon_l$ for each of the 12 eigenvalues in this illustrative example, where $\hat{\lambda}_1, \dots, \hat{\lambda}_4$ are the true eigenvalues marked by blue circles, and the rest are the spurious eigenvalues marked by red crosses.
Evidently, the absolute errors of $\hat{\lambda}_1, \dots, \hat{\lambda}_4$ are smaller by at least 10 orders of magnitude than the absolute errors of the spurious eigenvalues. Consequently, the true and spurious eigenvalues can be easily distinguished by visual inspection based on Fig. \ref{fig:abs_err_example}.

To implement this identification in an unsupervised fashion, $\varepsilon_l$ are clustered into two subsets, $S_1$ and $S_2$, using the $k$-means algorithm with $k=2$. The resulting clusters are $\{ \varepsilon_1, \dots, \varepsilon_4 \} \in S_1$ and $\{ \varepsilon_5, \dots, \varepsilon_{12} \} \in S_2$.
Next, the respective average absolute errors of $S_1$ and $S_2$ are computed, namely, $\langle \varepsilon \rangle _{S_1}$ and $ \langle \varepsilon \rangle _{S_2}$, according to \eqref{eq:avg_eig_error}. In this example, $\langle \varepsilon \rangle _{S_1} \approx 10^{-15}$ and $\langle \varepsilon \rangle _{S_2} \approx 10^{-3}$.
Therefore, subset $S_1$, corresponding to the smaller value $\langle \varepsilon \rangle _{S_1} < \langle \varepsilon \rangle _{S_2}$ is identified as the subset containing the $\varepsilon_l$ that are related to the true DMD eigenvalues.

Based on this identification, the original (sampled) signal is reconstructed using Eq. \eqref{eq:xk_analogous}. The accuracy of the reconstruction can be observed in Fig. \ref{fig:reconstructed}, which shows the original (black solid line) and the reconstructed (blue dashed line) signals.
We note that the markings of the blue circles and red crosses in Fig. \ref{fig:abs_err_example}, as well as the reconstructed signal in Fig. \ref{fig:reconstructed}, were obtained in a completely unsupervised and automatic manner using Algorithm \ref{alg}.

To conclude, we emphasize that the input of our algorithm is the samples (observations) of the system, where it operates in an \emph{equation-free} and \emph{unsupervised} fashion to, eventually, extract a reduced-order and optimal spectral representation of the system.

The code that reproduces the results in this section is openly available in the following \href{https://github.com/emilbronstein/spatiotemporal_coupling_in_delay_coordinates_DMD}{GitHub link}.


\section{Simulation results} \label{sec:results}


\subsection{Two-mode sine signal}\label{subsec:twomode}

Consider the two-mode sine signal in Eq. \eqref{eq:twomode}. We examine noise with four different amplitudes $\alpha$ that correspond to the signal-to-noise ratios (SNR) 15, 10, 5, and 0 dB.
We test the performance of the proposed method (detailed in section \ref{sec:proposed_method} and summarized in Algorithm \ref{alg}) for the augmentation numbers $s=10, 20, \dots, 360$, and compare its DMD components selection to the selection method based on maximal amplitudes.
For comparison purposes, we focus on the accuracy of the recovery of the oscillation frequencies of the signal obtained by the two methods.
For a fair comparison, we set the number of selected DMD components to 4 (2 oscillation frequencies that appear in 2 complex conjugate pairs correspond to 4 DMD eigenvalues). Note that, in general, our method does not require the number of DMD components as a prior, but infers it from the observations.

\begin{figure*}[t]
\includegraphics[width=0.7\textwidth]{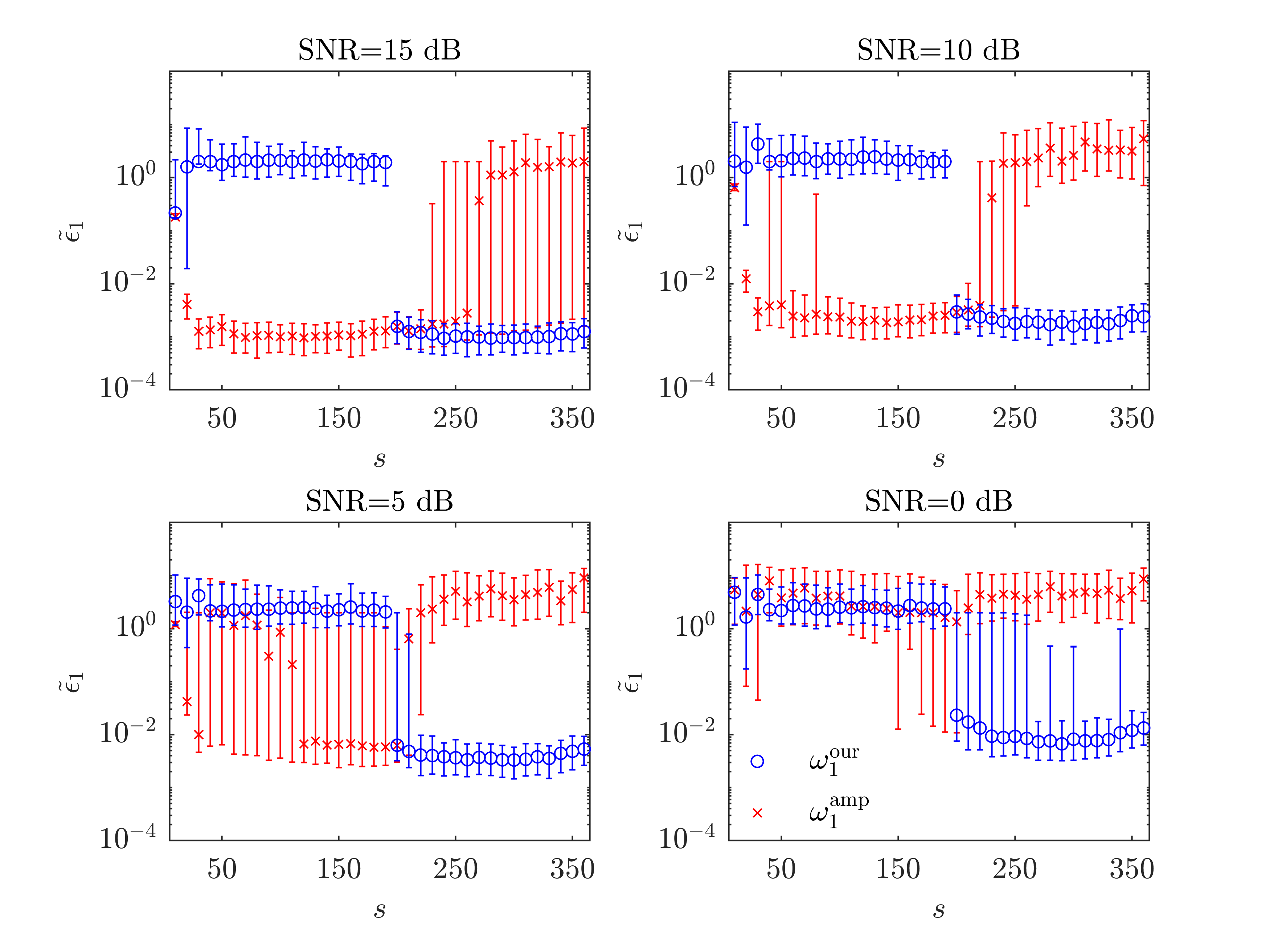}
\caption{\label{fig:e_two_mode_sine_w1} The medians of the absolute errors, $\Tilde{\epsilon}_1$, for augmentation numbers $s=10, 20, \dots, 360$ and SNR values of 15, 10, 5, 0 dB. The blue circles and red crosses denote $\Tilde{\epsilon}_1$ computed based on our and the maximal amplitudes methods, respectively. The whiskers denote the 25\textsuperscript{th} and 75\textsuperscript{th} percentiles of ${\epsilon}_1$.
}
\end{figure*}

\begin{figure*}[t]
\includegraphics[width=0.7\textwidth]{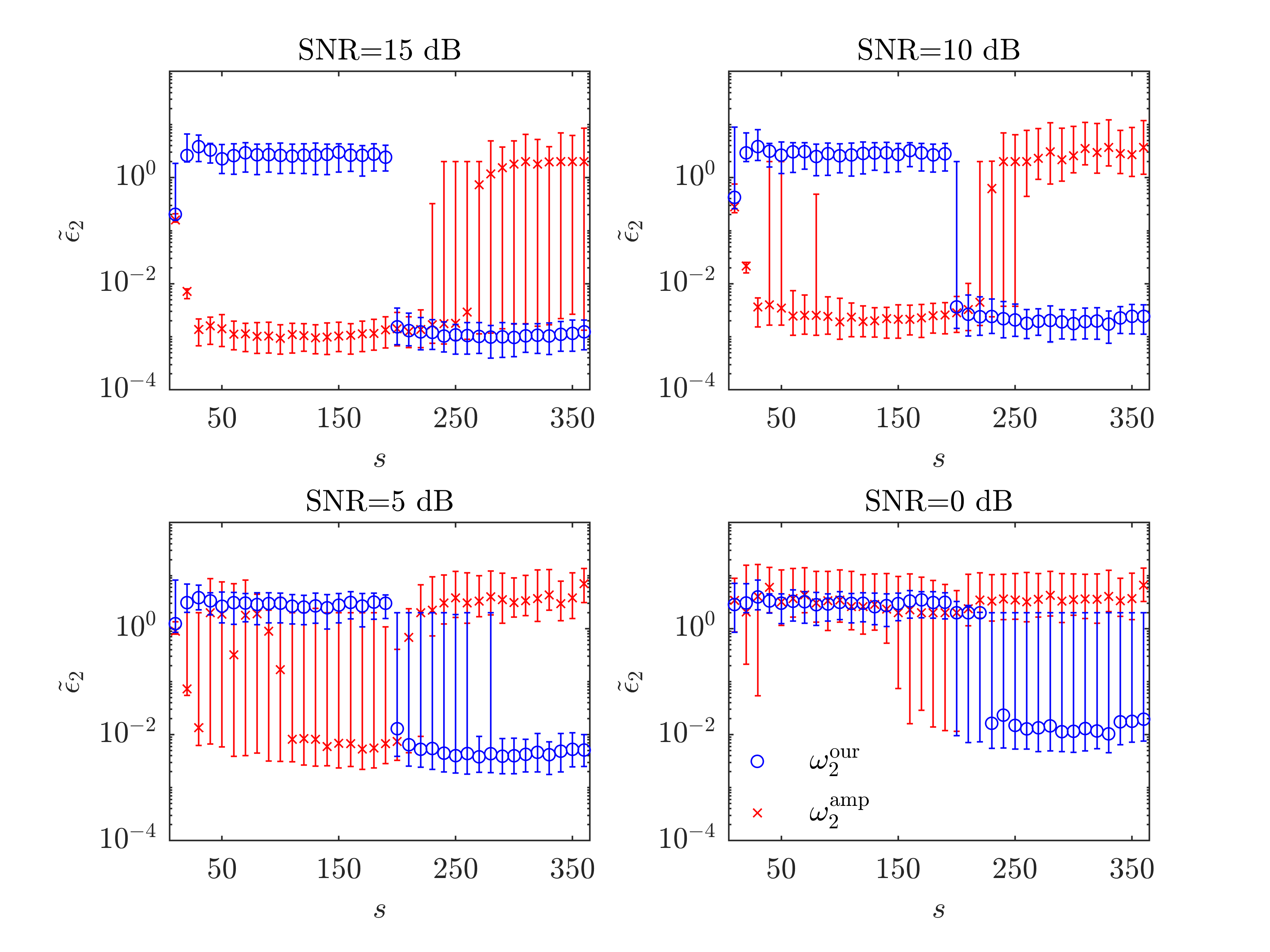}
\caption{\label{fig:e_two_mode_sine_w2}
The medians of the absolute errors, $\Tilde{\epsilon}_2$, for augmentation numbers $s=10, 20, \dots, 360$ and SNR values of 15, 10, 5, 0 dB. The blue circles and red crosses denote $\Tilde{\epsilon}_2$ computed based on our and the maximal amplitudes methods, respectively. The whiskers denote the 25\textsuperscript{th} and 75\textsuperscript{th} percentiles of ${\epsilon}_2$.
}
\end{figure*}

The recovery of the oscillation frequencies is evaluated by the absolute errors, $\epsilon_i$, between the system's true oscillation frequencies ($\omega_1^{\text{sys}}$ and $\omega_2^{\text{sys}}$) and the frequencies obtained from the DMD eigenvalues that were selected based on our and the maximal amplitudes methods, denoted by the superscripts \emph{our} and \emph{amp}, respectively. Namely, we calculate
\begin{equation}
    \epsilon_i^{\text{our}} = \left| \omega_i^{\text{sys}} - \omega_i^{\text{our}} \right|,
    \quad
    \epsilon_i^{\text{amp}} = \left| \omega_i^{\text{sys}} - \omega_i^{\text{amp}} \right|.
\end{equation}
We repeat the calculations of $\epsilon_i$ 500 times, each time with a different generated random noise, and compute the medians of the absolute errors, $\Tilde{\epsilon}_i$, as well as their 25\textsuperscript{th} and 75\textsuperscript{th} percentiles.
Fig. \ref{fig:e_two_mode_sine_w1} considers the first oscillation frequency, $\omega_1^{\text{sys}}=3$ rad/s, and presents $\Tilde{\epsilon}_1$ for different augmentation numbers $s$ at different SNR values, where the 25\textsuperscript{th} and 75\textsuperscript{th} percentiles of $\epsilon_1$ are denoted by the whiskers.
Fig. \ref{fig:e_two_mode_sine_w2} shows the same analysis for $\omega_2^{\text{sys}}=5$ rad/s.

Fig. \ref{fig:e_two_mode_sine_w1} and Fig. \ref{fig:e_two_mode_sine_w2} show that at the larger SNR values (15 and 10 dB), both $\omega_1^{\text{sys}}$ and $\omega_2^{\text{sys}}$ can be recovered with small errors based on both methods, given that $s$ is chosen according to the range indicated in \eqref{eq:s_range} for our method, and in \eqref{eq:s_range_small} for the maximal amplitudes method.

Yet, this is not the case for the smaller SNR values. For SNR value of 5 dB, the maximal amplitudes method is unable to recover $\omega_1^{\text{sys}}$ and $\omega_2^{\text{sys}}$ with reasonable errors; namely, the frequencies are obtained with either large error median or large error variance of $\epsilon_i$.
Contrarily, our method recovers these frequencies with small errors when the values of $s$ are set in  range \eqref{eq:s_range}.
For SNR value of 0 dB, the maximal amplitudes method leads to large error medians, whereas our method obtains small error medians when the values of $s$ are set in  range \eqref{eq:s_range}.


\subsection{Quasiperiodic signal}\label{subsec:quasiperiodic}

Consider the following noisy quasiperiodic signal

\begin{equation}
x(t) = \sin \left(\sqrt{10}t \right) \sin(t) + \alpha n(t),
\end{equation}
which can be recast as
\label{eq:quasiperodic}
\begin{equation}
x(t) = \frac{1}{2} \cos \left((\sqrt{10}-1)t \right) -  \frac{1}{2} \cos \left((\sqrt{10}+1)t \right) + \alpha n(t),
\label{eq:quasiperodic_alternative}
\end{equation}
where $n(t)$ is a white standard Gaussian noise with varying amplitudes $\alpha$ that correspond to SNR value of 15, 10, 5, and 0 dB.

\begin{figure*}[t]
\includegraphics[width=0.7\textwidth]{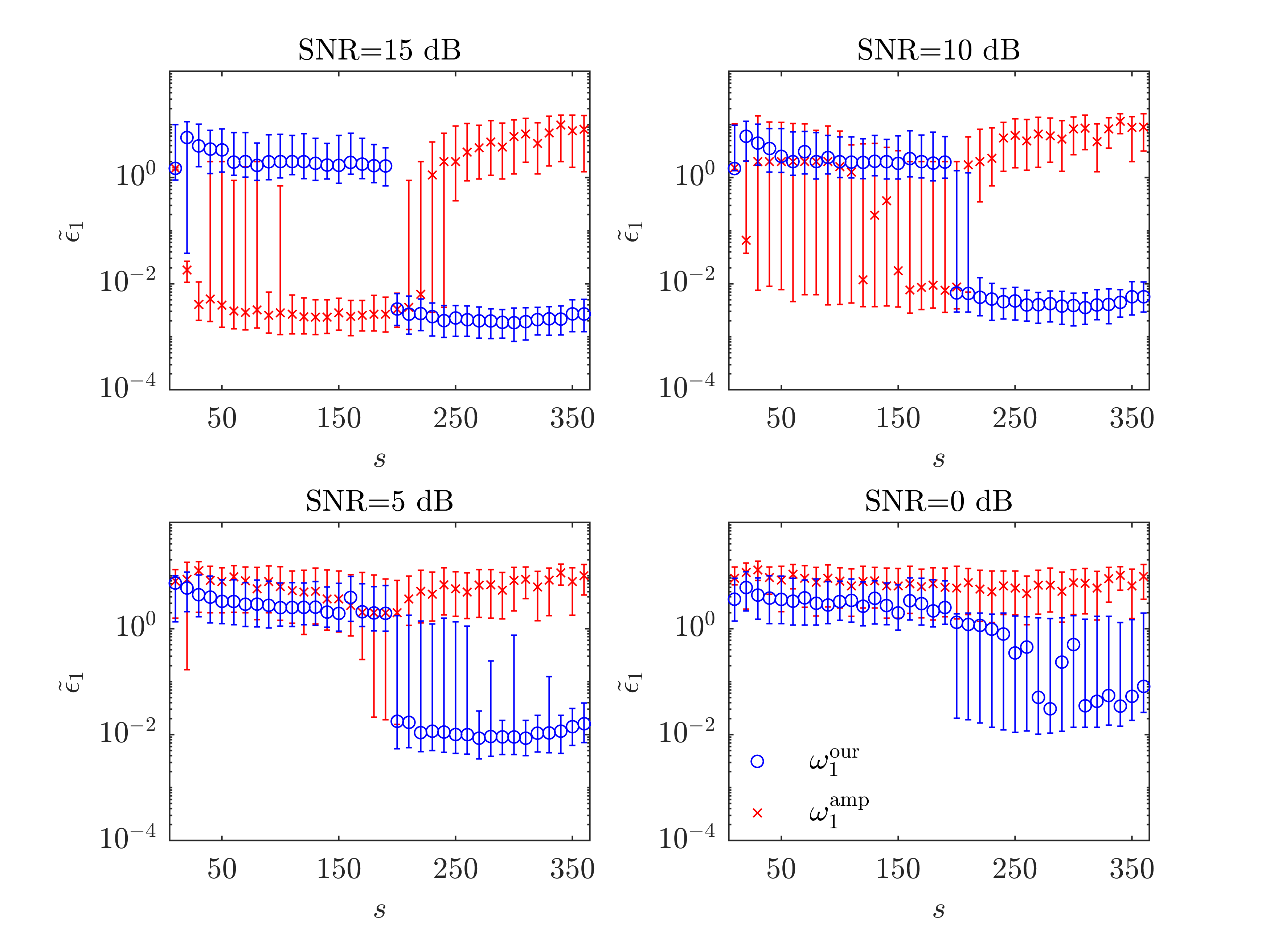}
\caption{\label{fig:quasiperiodic_e_w1}
Same as Fig. \ref{fig:e_two_mode_sine_w1}, but for $\omega_1$ of the quasiperiodic signal.
}
\end{figure*}

\begin{figure*}[t]
\includegraphics[width=0.7\textwidth]{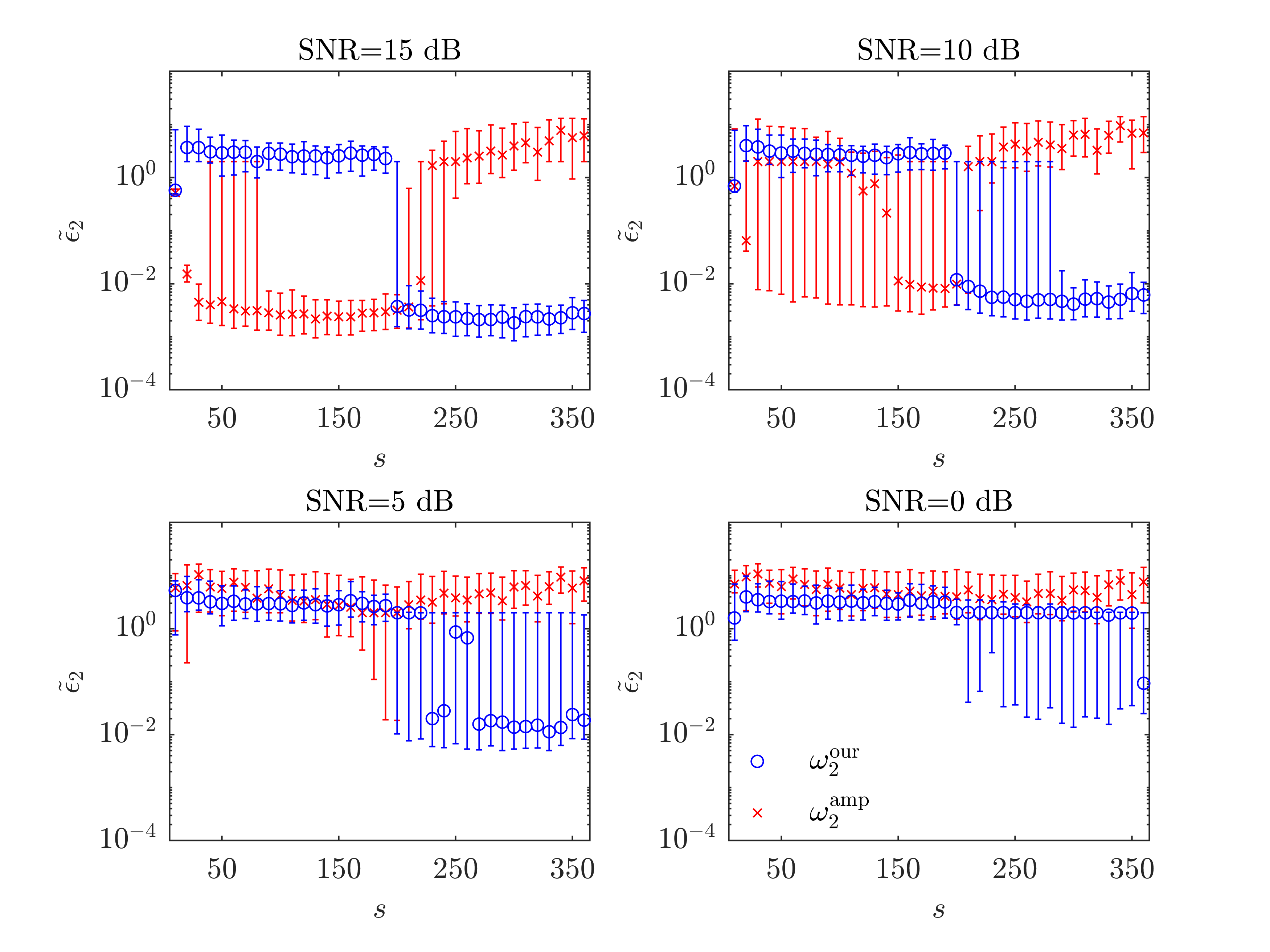}
\caption{\label{fig:quasiperiodic_e_w2}
Same as Fig. \ref{fig:e_two_mode_sine_w2}, but for $\omega_2$ of the quasiperiodic signal.
}
\end{figure*}

Similarly to the two-mode sine signal simulation (section \ref{subsec:twomode}), we test the capability to uncover the signal's irrational frequencies in \eqref{eq:quasiperodic_alternative} using our method and compare it with the maximal amplitudes method. Fig. \ref{fig:quasiperiodic_e_w1} and \ref{fig:quasiperiodic_e_w2} are the same as Fig. \ref{fig:e_two_mode_sine_w1} and \ref{fig:e_two_mode_sine_w2} but with respect to $\omega_1^{\text{sys}} = ( \sqrt{10}-1 )$ rad/s and $\omega_2^{\text{sys}} = ( \sqrt{10}+1 )$ rad/s, respectively.

For both $\omega_1^{\text{sys}}$ and $\omega_2^{\text{sys}}$, the maximal amplitudes method performs well only at SNR value of 15 dB. However, at the other three examined SNR values, the system's frequencies are uncovered with a large error variance (10 dB) or a large error median (5 and 0 dB).
Contrarily, our method uncovers $\omega_1^{\text{sys}}$ with small errors at SNR=15, 10, and 5 dB, and with small error median but with large variance at 0 dB. Also, using our method, $\omega_2^{\text{sys}}$ is uncovered with small errors at SNR=15 and 10 dB, and small error median but with large variance at 5 dB. At SNR = 0 dB, both methods fail to uncover $\omega_2^{\text{sys}}$.

\subsection{Multi degrees of freedom oscillator}\label{subsec:oscillator}

\begin{figure*}

\begin{tikzpicture}[black!75,thick]
 
\fill [pattern = north west lines] (0,0) rectangle (10,0.5);
\draw[thick] (0,0.5) -- (10,0.5);

\fill [pattern = north west lines] (-0.5,0) rectangle (0,2);
\draw[thick] (0,0.5) -- (0,2);

\draw
[
    decoration={
        coil,
        aspect=0.3, 
        segment length=1.2mm, 
        amplitude=2mm, 
        pre length=3mm,
        post length=3mm},
    decorate
] (0,1.25) -- (2,1.25) 
    node[midway,above=0.25cm,black]{$k_1$};

\node[draw,
    fill=white!60,
    minimum width=1cm,
    minimum height=1cm,
    anchor=north,
    label=center:$m_1$] at (2.5,1.75) {};
    \draw [thick, fill=black] (2.25,0.62) circle (0.1) (2.75,0.62) circle (0.1);

\draw[very thick,
    black,
    |-latex] (2.5,2) -- (3.5,2)
    node[midway,above]{$x_1(t)$};

\draw
[
    decoration={
        coil,
        aspect=0.3, 
        segment length=1.2mm, 
        amplitude=2mm, 
        pre length=3mm,
        post length=3mm},
    decorate
] (3,1.25) -- (5,1.25) 
    node[midway,above=0.25cm,black]{$k_2$}; 

\node[draw,
    fill=white!60,
    minimum width=1cm,
    minimum height=1cm,
    anchor=north,
    label=center:$m_2$] at (5.5,1.75) {};
    \draw [thick, fill=black] (5.25,0.62) circle (0.1) (5.75,0.62) circle (0.1);

\draw[very thick,
    black,
    |-latex] (5.5,2) -- (6.5,2)
    node[midway,above]{$x_2(t)$};

\draw
[
    decoration={
        coil,
        aspect=0.3, 
        segment length=1.2mm, 
        amplitude=2mm, 
        pre length=3mm,
        post length=3mm},
    decorate
] (6,1.25) -- (8,1.25) 
    node[midway,above=0.25cm,black]{$k_3$};  
    
\node[draw,
    fill=white!60,
    minimum width=1cm,
    minimum height=1cm,
    anchor=north,
    label=center:$m_3$] at (8.5,1.75) {};
    \draw [thick, fill=black] (8.25,0.62) circle (0.1) (8.75,0.62) circle (0.1);

\draw[very thick,
    black,
    |-latex] (8.5,2) -- (9.5,2)
    node[midway,above]{$x_3(t)$};

\end{tikzpicture}

\caption{\label{fig:undamped_oscillator}
Three degrees of freedom oscillator, comprised of three masses $m_1=1$ kg, $m_2=2$ kg, $m_3=3$ kg, connected via three springs of constants $k_1=k_2=50$ N/m, $k_3 = 75$ N/m. The system is subjected to the initial conditions $x_1(0) = 1$ m, $x_2(0) = 2$ m, $x_3(0) = 3$ m, $\dot{x}_1(0)=\dot{x}_2(0)=\dot{x}_3(0)=0$, giving rise to its oscillations.
}

\end{figure*}
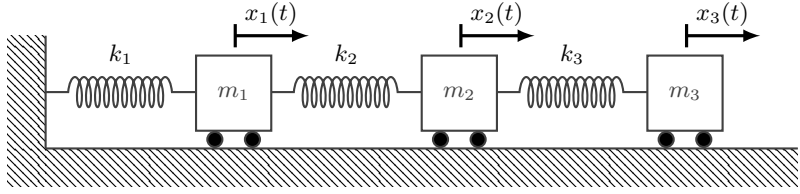

Consider an oscillator comprised of $N$ masses $m_1, \dots, m_N$ connected by springs of constants $k_1, \dots, k_N$. The first mass $m_1$ is connected to a wall via the first spring $k_1$, as depicted in Fig. \ref{fig:undamped_oscillator} for $N=3$. The equations of motion of this system are given by 

\begin{equation*}
    \bm{M} \ddot{\bm{x}}(t)
    +
    \bm{K} \bm{x}(t)
    =
    \bm{0},
\end{equation*}
where
\begin{equation*}
\begin{split}
    \bm{M}
    & =
    \begin{bmatrix}
    m_1 & 0 & \cdots & 0 \\
    0 & m_2 & \cdots & 0 \\
    \vdots & \vdots & \ddots & \vdots \\
    0 & 0 & \cdots & m_N
    \end{bmatrix},
    \\
    \bm{K}
    & =
    \begin{bmatrix}
    k_1+k_2 & -k_2 & 0 & \cdots & 0 \\
    -k_2 & k_2+k_3 & -k_3 & \cdots & 0 \\
    \vdots & \ddots & \ddots & \ddots & \vdots \\
    0 & 0 & 0 & -k_N & k_N
    \end{bmatrix},
\end{split}
\end{equation*}
and $\bm{x} = [x_1, x_2, \dots, x_N]^T$ are the positions of the masses.

To excite motion, the system is perturbed from its rest state by moving each mass to some initial position, and then releasing them all at once. Namely, the system is subjected to the initial conditions $\bm{x}(t=0) = [x_1^0, x_2^0, \dots, x_N^0]^T$, as well as to zero initial velocities of all masses ($\dot{\bm{x}}(t=0) = \bm{0}$).

Here, we consider the 3 DOF oscillator in Fig. \ref{fig:undamped_oscillator} ($N=3$) with the masses $m_1=1$ kg, $m_2=2$ kg, $m_3=3$ kg, spring constants $k_1=k_2=50$ N/m, $k_3 = 75$ N/m, initial positions $x_1(0) = 1$ m, $x_2(0) = 2$ m, $x_3(0) = 3$ m, and zero initial velocities $\dot{x}_1(0)=\dot{x}_2(0)=\dot{x}_3(0)=0$. In this case, the system's natural frequencies are \cite{meirovitch2001fundamentals}
\begin{equation}
\begin{split}
    \omega_1^{\text{sys}} & = 2.0626 \ \text{rad/s},
    \\
    \omega_2^{\text{sys}} & = 7.6961 \ \text{rad/s},
    \\
    \omega_3^{\text{sys}} & = 11.1363 \ \text{rad/s}.
\end{split}
    \label{3dof_w}
\end{equation}

\begin{figure*}[t]
\includegraphics[width=0.7\textwidth]{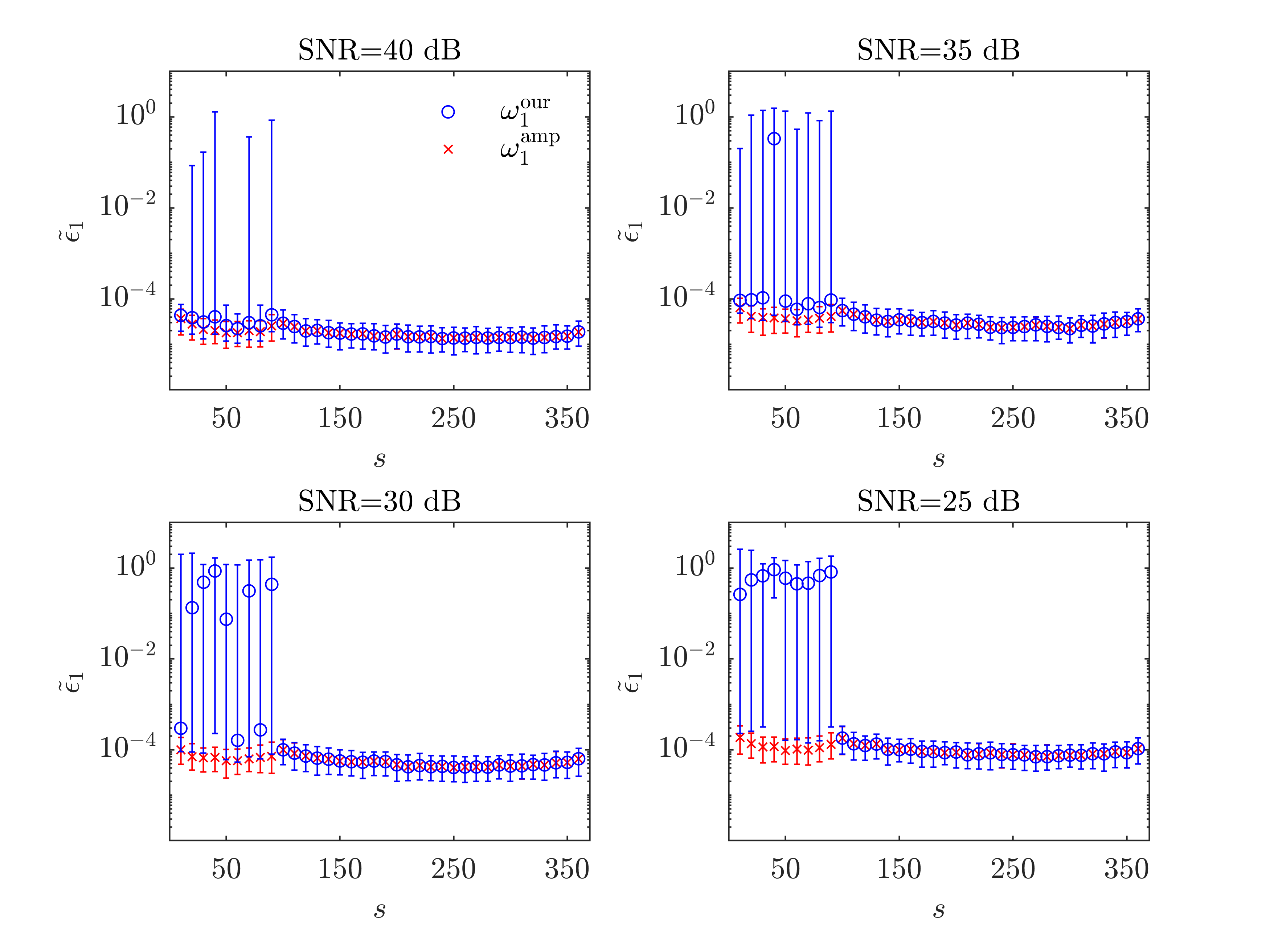}
\caption{\label{fig:undamped_oscillator_e_w1}
Same as Fig. \ref{fig:e_two_mode_sine_w1}, but for $\omega_1$ of the 3 DOF oscillator.
}
\end{figure*}

\begin{figure*}[t]
\includegraphics[width=0.7\textwidth]{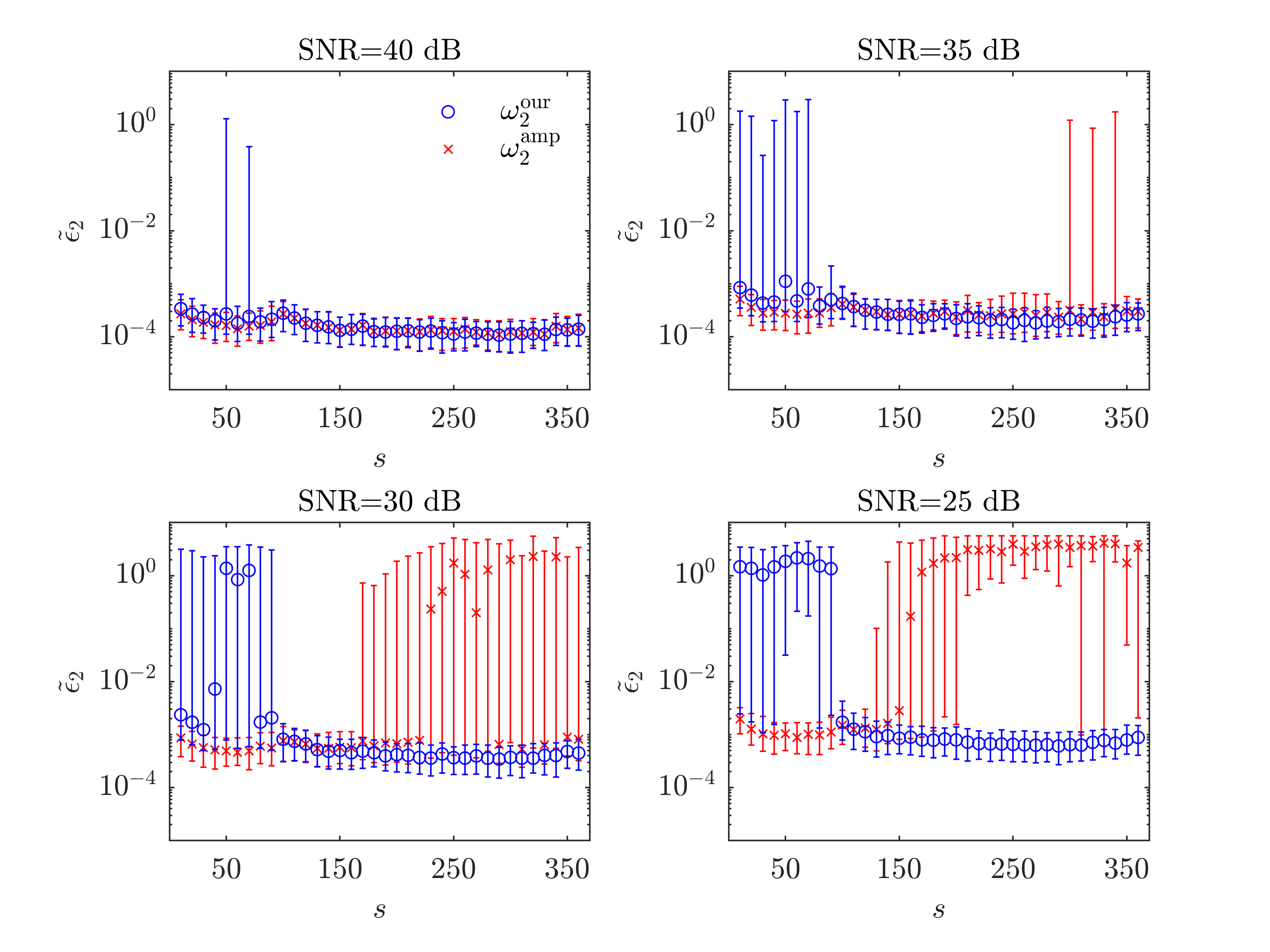}
\caption{\label{fig:undamped_oscillator_e_w2}
Same as Fig. \ref{fig:e_two_mode_sine_w1}, but for $\omega_2$ of the 3 DOF oscillator.
}
\end{figure*}

\begin{figure*}[t]
\includegraphics[width=0.7\textwidth]{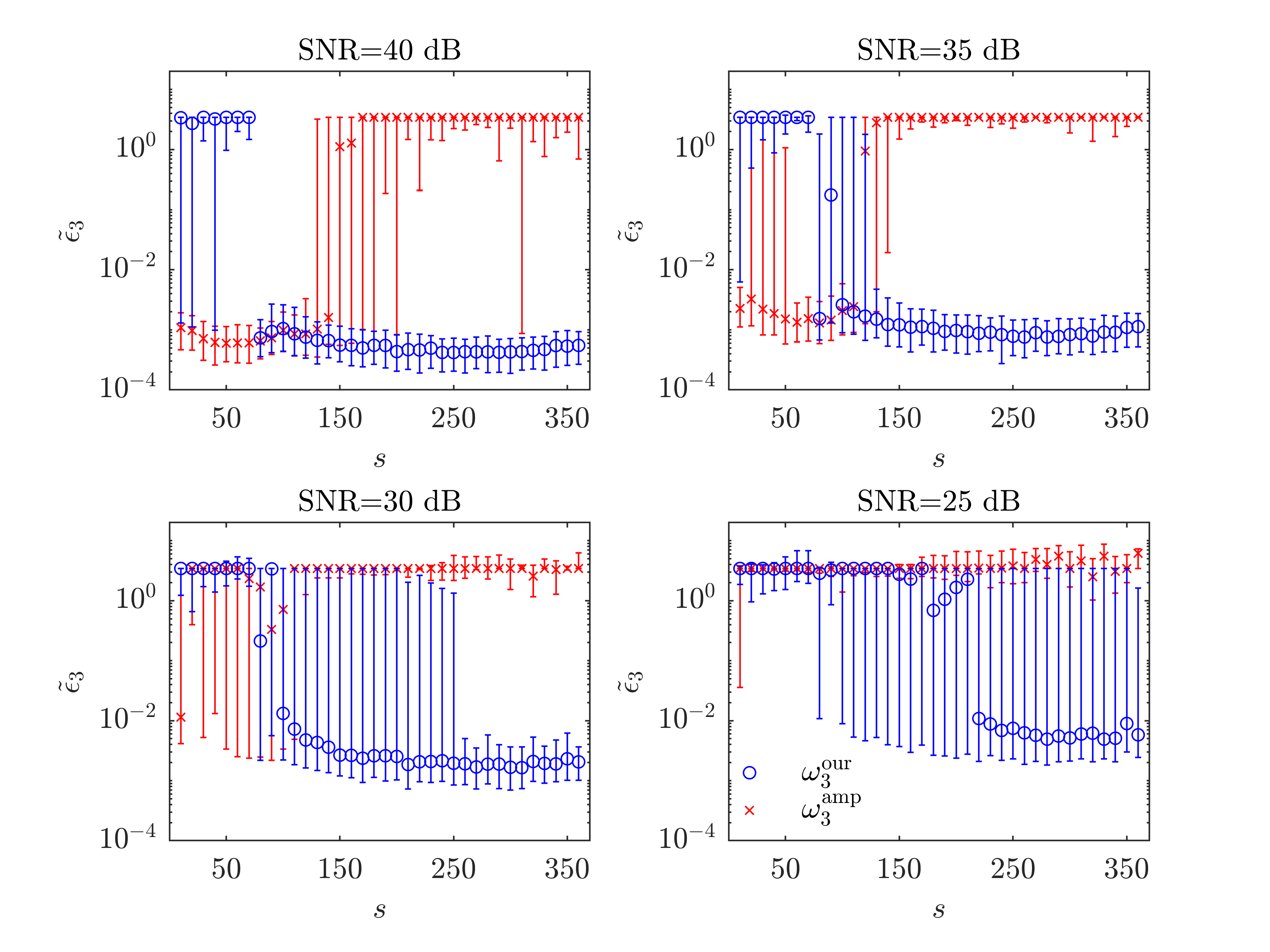}
\caption{\label{fig:undamped_oscillator_e_w3}
Same as Fig. \ref{fig:e_two_mode_sine_w1}, but for $\omega_3$ of the 3 DOF oscillator.
}
\end{figure*}

Random noise drawn from a normal distribution with amplitudes that correspond to SNR value of 40, 35, 30, and 25 dB is added to the observations. Note that in this case, the SNR is computed with respect only to the first mass. Then, same as in Fig. \ref{fig:e_two_mode_sine_w1}, the accuracy of the recovery of frequencies \eqref{3dof_w} by our and the maximal amplitudes methods is tested and presented in Fig. \ref{fig:undamped_oscillator_e_w1}, Fig. \ref{fig:undamped_oscillator_e_w2}, and Fig. \ref{fig:undamped_oscillator_e_w3}.

At all the examined SNR values, $\omega_1^{\text{sys}}$ and $\omega_2^{\text{sys}}$ can be recovered with small errors using both methods. These results are similar for $\omega_3^{\text{sys}}$ at SNR values of 40 and 35 dB.
However, the maximal amplitudes method fails to uncover $\omega_3^{\text{sys}}$ with reasonable errors at SNR value of 30 dB, while our method successfully performs this task. At SNR value of 25 dB, the maximal amplitudes method is unable to uncover $\omega_3^{\text{sys}}$ as well, whereas our method uncovers this frequency with small error median and large error variance.


\subsection{Damped oscillator}\label{subsec:damped_oscillator}

\begin{figure}

\begin{tikzpicture}[black!75,thick]

\tikzstyle{damper}=[thick,decoration={markings,  
  mark connection node=dmp,
  mark=at position 0.5 with 
  {
    \node (dmp) [thick,inner sep=0pt,transform shape,rotate=-90,minimum
width=15pt,minimum height=3pt,draw=none] {};
    \draw [thick] ($(dmp.north east)+(2pt,0)$) -- (dmp.south east) -- (dmp.south
west) -- ($(dmp.north west)+(2pt,0)$);
    \draw [thick] ($(dmp.north)+(0,-5pt)$) -- ($(dmp.north)+(0,5pt)$);
  }
}, decorate]
 
\fill [pattern = north west lines] (0,0) rectangle (5,0.5);
\draw[thick] (0,0.5) -- (5,0.5);

\fill [pattern = north west lines] (-0.5,0) rectangle (0,3);
\draw[thick] (0,0.5) -- (0,3);

\draw
[
    decoration={
        coil,
        aspect=0.3, 
        segment length=1.2mm, 
        amplitude=2mm, 
        pre length=3mm,
        post length=3mm},
    decorate
] (0,2) -- (2,2) 
    node[midway,above=0.25cm,black]{$k$}; 

\draw [damper] ($(0,1.2)$) -- ($(2,1.2) $)
node[midway,below=0.25cm,black]{$c$}; 

\node[draw,
    fill=white!60,
    minimum width=1.5cm,
    minimum height=1.5cm,
    anchor=north,
    label=center:$m$] at (2.75,2.25) {};
\draw [thick, fill=black] (2.25,0.62) circle (0.1) (3.25,0.62) circle (0.1) (2.75,0.62) circle (0.1);

\draw[very thick,
    black,
    |-latex] (2.75,2.5) -- (3.75,2.5)
    node[midway,above]{$x(t)$};

\end{tikzpicture}

\caption{\label{fig:single_damped_mass}
A damped oscillator comprised of mass $m=1$ kg, which is connected to a wall via a spring of stiffness  $k=49$ N/m, and a damper with a damping coefficient $c= 0.07$ Ns/m. The mass is subjected to the initial conditions $x(0) = 1$ m, $\dot{x}(0)=2$ m/s, giving rise to its oscillations.
}

\end{figure}
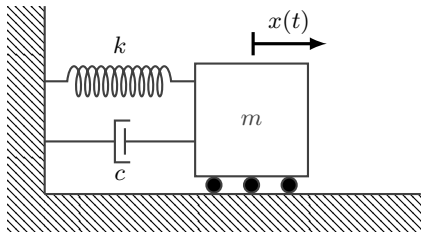

Consider a damped oscillator with a mass $m = 1$ kg, connected to a wall via a spring with constant $k=49$ N/m and a damper with damping coefficient $c= 0.07$ Ns/m, as illustrated in Fig. \ref{fig:single_damped_mass}. These system parameters correspond to the natural frequency $\omega_n = 7$ rad/s, damping ratio $\zeta=0.005$, and the frequency of damped vibrations $\omega_d = \omega_n \sqrt{1- \zeta^2} = 6.9999$ rad/s. The mass is subjected to the initial conditions $x(0) = x_0 = 1$ m, $\dot{x}(0) = v_0 = 2$ m/s, giving rise to oscillations, which induce the following time response \cite{meirovitch2001fundamentals}

\begin{equation}
    x(t) = \exp(-\zeta \omega_n t) \left(\frac{v_0 + \zeta \omega_n x_0}{\omega_d} \sin (\omega_d t)  +  x_0 \cos(\omega_d t) \right).
    \label{eq:damped_oscillator}
\end{equation}

\begin{figure*}[t]
\includegraphics[width=0.7\textwidth]{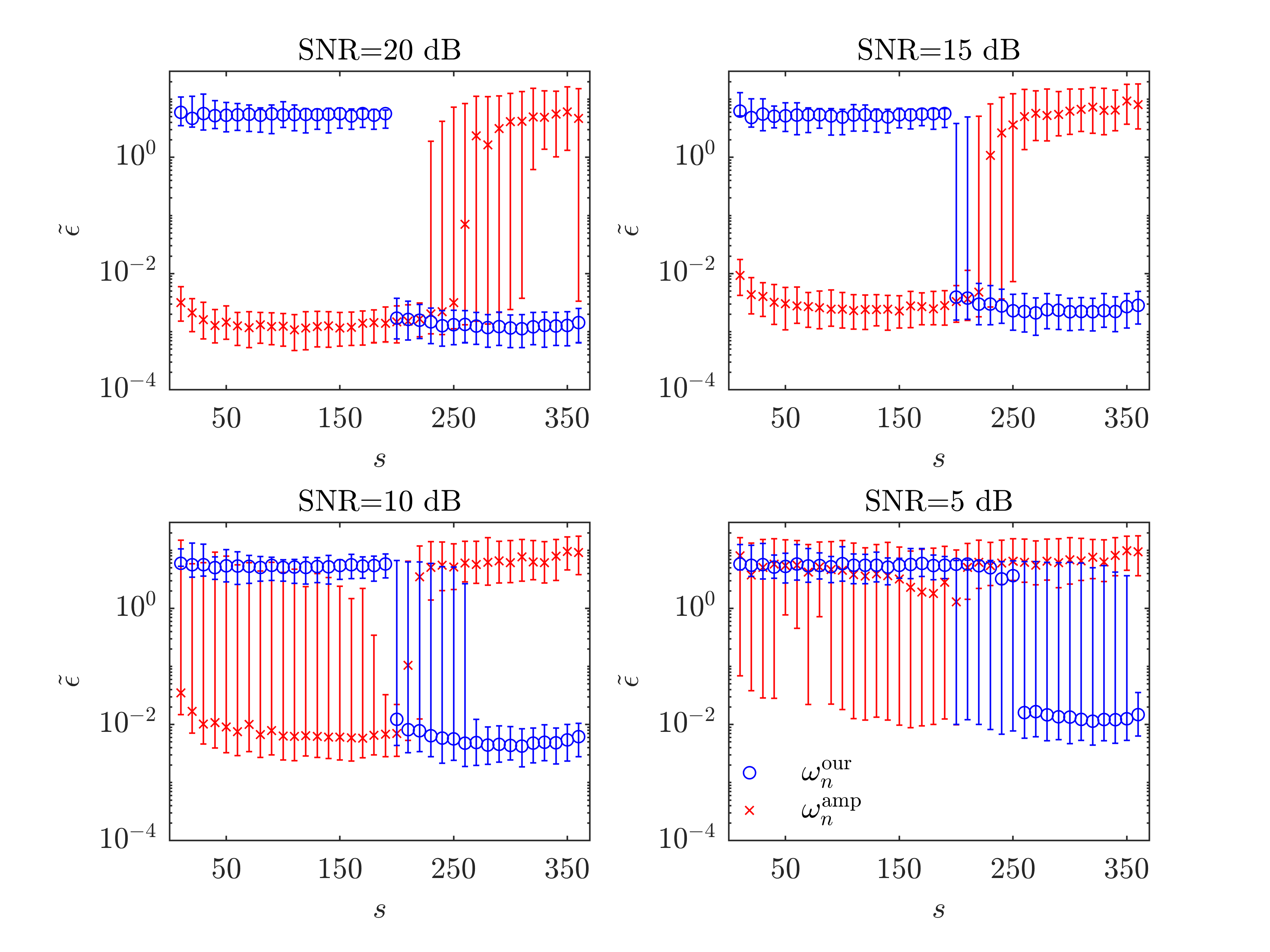}
\caption{\label{fig:single_damped_mass_e_w}
Same as Fig. \ref{fig:e_two_mode_sine_w1}, but for $\omega_n$ of the damped oscillator.
}
\end{figure*}

\begin{figure*}[t]
\includegraphics[width=0.7\textwidth]{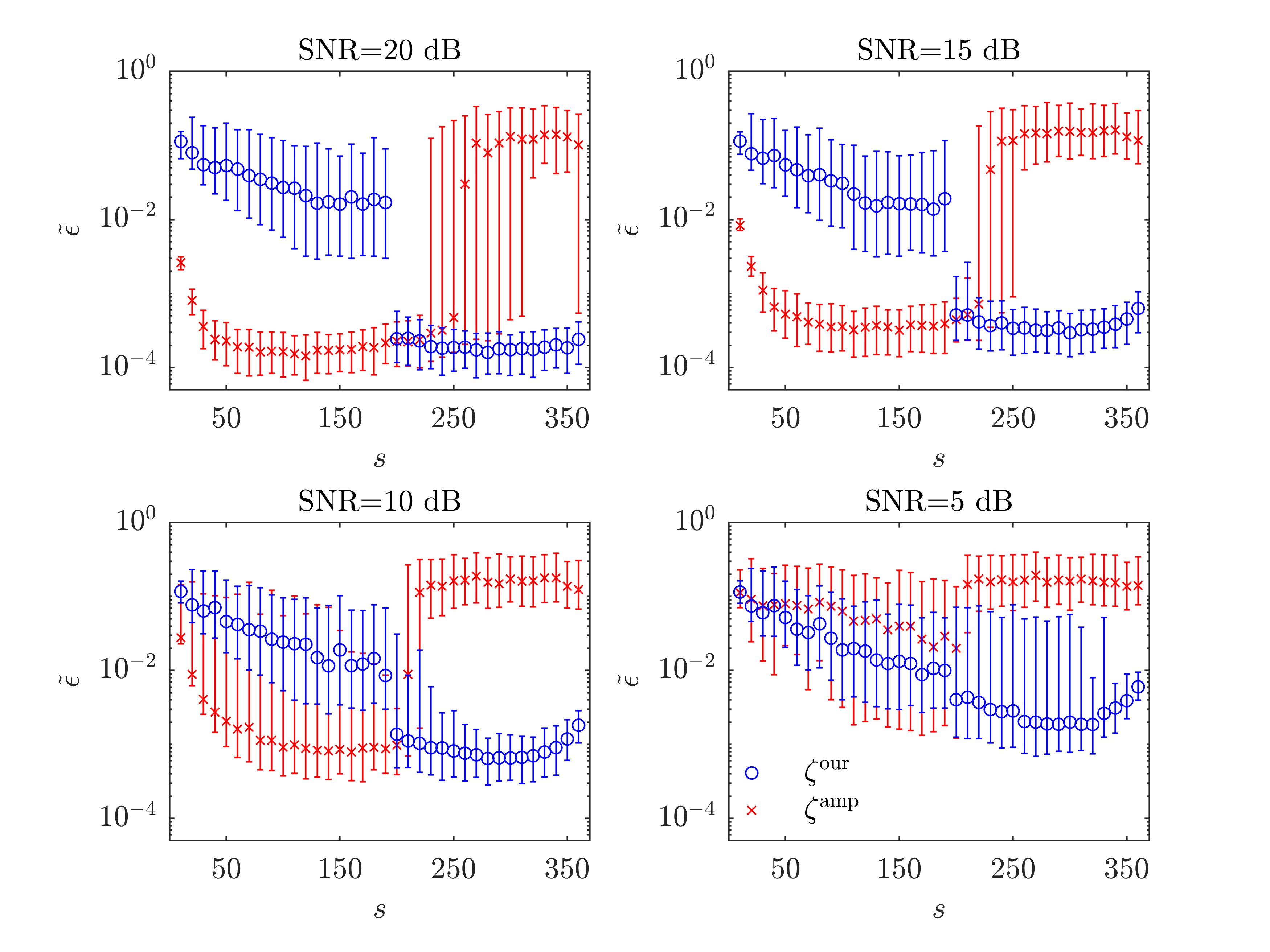}
\caption{\label{fig:single_damped_mass_e_z}
Same as Fig. \ref{fig:e_two_mode_sine_w1}, but for $\zeta$ of the damped oscillator.
}
\end{figure*}

Random noise drawn from a normal distribution with amplitudes that correspond to SNR values of  20, 15, 10, and 5 dB is added to \eqref{eq:damped_oscillator}. The recovery accuracy of $\omega_n$ and $\zeta$ obtained by our method and the maximal amplitudes method is presented in Fig. \ref{fig:single_damped_mass_e_w} and Fig. \ref{fig:single_damped_mass_e_z}, respectively.

At SNR values of 20 and 15 dB, $\omega_n$ and $\zeta$ are recovered with small errors using both methods. At SNR values of 10 and 5 dB, the maximal amplitudes method is unable to recover both $\omega_n$ and $\zeta$ with reasonable errors. Conversely, at these noise levels, our method yields a recovery with errors at the order of $10^{-2}$ for $\omega_n$ and $10^{-3}$ for $\zeta$.


\section{Conclusions}

Delay-coordinates DMD is widely-used for data-driven analysis of dynamical systems based on observations in a broad range of fields.
In this paper, we investigated two key questions concerning the inherent redundancy that arises from the utility of delay-coordinates DMD: the excess of dynamical components (i.e., spurious components) and the excess of dimensionality (i.e., coordinates).
We showed that delay-coordinates DMD induces a particular structure on the augmented DMD components, consisting of a spatiotemporal coupling.
At first glance, this coupling seems to counter the core idea underlying DMD, which facilitates a representation of the system that decouples temporal and spatial patterns.
Yet, a deeper look allowed us not only to mitigate this coupling, but also to exploit it.
Based on the spatiotemporal coupling we presented, we proposed a method for constructing a compact and improved reduced-order DMD representation.
Specifically, we showed how to identify and select the informative (true) DMD components, thereby addressing the excess of dynamical components. This identification is based on the induced temporal associations within each augmented mode, which allowed us to address the redundancy in dimensionality. 
We tested the proposed method on four dynamical systems corrupted with noise and compared the performance to the prevalent method, which is based on the maximal amplitudes of the DMD modes. The results demonstrate the advantages of the proposed method, especially at the presence of high levels of noise.



%
%

%

\begin{acknowledgments}
This research was supported by the Pazy Foundation (grant No. 78-2018) and by the Israel Science Foundation (grants
No. 1309/18 and 1600/22).
E.B. is grateful to the Azrieli Foundation for the award of an Azrieli Fellowship.
R.T. acknowledges the support of the Schmidt Career Advancement Chair in AI.
\end{acknowledgments}

\bibliography{Bibliography.bib}

\end{document}